\documentclass[a4paper,11pt]{article}
\usepackage[top=2.5cm,bottom=2.5cm,left=2.2cm,right=2.2cm]{geometry}
\usepackage{amsfonts}
\usepackage{mathrsfs,amscd,amssymb,amsthm,amsmath,bm,graphicx,psfrag,subfigure,url,mathtools}
\usepackage{pict2e}
\usepackage{ulem}
\usepackage{stfloats}
\usepackage{psfrag,amsmath}
\usepackage{tikz}
\usepackage{indentfirst}
\usepackage{hyperref}
\usepackage{bookmark}
\usepackage{enumerate}
\usepackage{latexsym,euscript,epic,eepic,color}
\usepackage{multirow}
\usepackage{multicol}
\usepackage{longtable}
\usepackage{adjustbox}
\usepackage[all]{xy}
\usepackage{setspace}
\allowdisplaybreaks
\usepackage{authblk}
\usepackage{amsfonts}
\usepackage{amsmath}
\usepackage{mathrsfs,color}
\usepackage{epstopdf}
\usepackage{graphicx}
\usepackage{tikz}
\usetikzlibrary{positioning}
\usetikzlibrary{fadings}
\usetikzlibrary{patterns}
\usetikzlibrary{calc}
\usetikzlibrary{shadings}
\pgfsetlayers{background,main,foreground}

\makeatletter

\renewcommand{\@seccntformat}[1]{{\csname the#1\endcsname}{\normalsize .}\hspace{.5em}}
\makeatother

\usepackage{ifpdf}
\usepackage{indentfirst}

\def \[{\begin{equation}}
\def \]{\end{equation}}

\def\b{\color{blue}}

\newtheorem{thm}{Theorem}[section]

\newtheorem{claim}{Claim}
\newtheorem{lem}[thm]{Lemma}
\newtheorem{cor}[thm]{Corollary}

\newtheorem{pb}{Problem}

\newtheorem{conj}[thm]{Conjecture}

\newenvironment{wst}
{\setlength{\leftmargini}{1.5\parindent}
 \begin{itemize}
 \setlength{\itemsep}{-1.1mm}}
{\end{itemize}}
\begin{document}
\baselineskip=0.23in
\title{\bf Spectral extremal results on the $A_\alpha$-spectral radius of graphs without $K_{a,b}$-minor\thanks{S.L. is financially supported by the National Natural Science Foundation of China (Grant Nos. 12171190, 11671164) and the Special Fund for Basic Scientific Research of Central Colleges (Grant No. CCNU24JC005)\\[1pt] \hspace*{5mm}{\it Email addresses}: xyleiyuki@aliyun.com (X.Y. Lei),\  lscmath@ccnu.edu.cn (S.C. Li)}}
\author [1]{Xingyu Lei}
\author [,1,2]{{Shuchao Li}\thanks{Responding author}}
\affil[1]{School of Mathematics and Statistics, and Hubei Key Lab--Math. Sci.,\linebreak Central China Normal University, Wuhan 430079, China}
\affil[2]{Key Laboratory of Nonlinear Analysis \& Applications (Ministry of Education),\linebreak Central China Normal University, Wuhan 430079, China}
\date{\today}
\maketitle
\begin{abstract}
An important theorem about the spectral Tur\'an problem of $K_{a,b}$ was largely developed in separate papers. Recently it was completely resolved by Zhai and Lin [J. Comb. Theory, Ser. B 157 (2022) 184-215], which also confirms a conjecture proposed by Tait [J. Comb. Theory, Ser. A 166 (2019) 42-58]. Here, the prior work is fully stated, and then generalized with a self-contained proof. The more complete result is then used to better understand the relationship between the $A_{\alpha}$-spectral radius and the structure of the corresponding extremal $K_{a,b}$-minor free graph.
\vskip 0.2cm

\end{abstract}

\section{\normalsize Introduction}
We focus exclusively on simple, undirected graphs in this paper. Denote a graph as $G=(V(G), E(G))$, consisting of a {\it vertex set} $V(G)$ and an {\it edge set} $E(G)$. The number of vertices, and the number of edges, referred to as $n(G)$ and $e(G)$, are respectively known as the {\it order} and {\it size}  of the graph $G$. Unless mentioned otherwise, we will adhere to the standard notation and terminology commonly used in the field (see \cite{0001,0009}).

Let $A(G)$ and $D(G)$ be the adjacency matrix and degree diagonal matrix of $G$, respectively. The well-known \textit{signless Laplacian matrix} $Q(G)$ of $G$ satisfies $Q(G)=D(G)+A(G)$. The $A_\alpha$-\textit{matrix} of $G$ was introduced by Nikiforov \cite{Niki1} in 2017:
$$A_\alpha(G)=\alpha D(G)+(1-\alpha)A(G),\text{ where } \alpha\in[0,\,1].$$
The {\it characteristic polynomial} of $A_\alpha(G)$ is denoted by $\phi_{G}(x)$. The eigenvalues associated with $A_{\alpha}$ are referred to as the {\it $A_{\alpha}$-eigenvalues} of the graph $G$. The largest absolute value among these eigenvalues, represented as $\lambda_\alpha(G)$, is referred to as the {\it $A_{\alpha}$-spectral radius} of $G$. The study on the $A_{\alpha}$-spectral properties attracts more and more researchers' attention. For further developments on this subject, we direct readers to \cite{Li-Sun-1,Li-Sun-2,Li-Sun-3,Li-Wang,Li-wei,Lin,Niki2}.

A \textit{minor} of a graph $G$ refers to any graph that can be derived from $G$ by sequentially removing vertices and edges, as well as by contracting edges. Let $G$ and $H$ be two graphs. If $G$ does not include $H$ as a subgraph (resp. minor), then we call $G$ is \textit{$H$-free} (resp. \textit{$H$-minor free}).

The Tur\'an type problem is one of the central problems in graph theory. Similar to the Tur\'an type problem, Nikiforov introduced a spectral variation of this concept: 
What is the maximum $A_0$-spectral radius that an $H$-free graph with $n$ vertices may attains for a given graph $H$? This problem was studied extensively; see $K_a$ \cite{Boll,Wiff}, $K_{a,b}$ \cite{Baba,NikiK}, $P_k$ \cite{NikiP}, $C_4$ \cite{NikiC,Zhai-Wang}, $W_{2k+1}$ \cite{SB-1}, $F_k$ \cite{SB-2}, $\bigcup_{i=1}^kS_{r_i}$ \cite{Z1}, $\bigcup_{i=1}^kP_{r_i}$ \cite{Z2}, consecutive cycles \cite{G1}, consecutive short odd cycles \cite{L1,Lin1}, disjoint cliques \cite{N1}, disjoint cycles \cite{Fang} and some nice surveys \cite{Chen-Zhang,Li-Liu,NikiA}.

It is obvious that $A_0(G)=A(G),\,A_{\frac{1}{2}}(G)=\frac{1}{2}Q(G)$ and $A_1(G)=D(G).$ In 2019, Tait \cite{Tait} put forward the following conjecture.
\begin{conj}[\cite{Tait}]\label{conjmain}
Given $b \geqslant a\geqslant 2$ and $n-1=kb+t\ (0\leqslant t\leqslant b-1)$, then $K_{a-1}\vee (k K_b\cup K_t)$ attains the maximum $A_0$-spectral radius among all $K_{a,b}$-minor free graphs with $n$ vertices for large enough $n$.
\end{conj}
In the past several years, Conjecture \ref{conjmain} was confirmed for $a + b = 4$ \cite{NikiC,Zhai-Wang}, $a + b = 5$ \cite{Niki-2t}, $a + b = 6$ \cite{Wang}, and $t = 0$ \cite{Tait}. Recently, Zhai and Lin \cite{Zhai-Lin} completely resolved the rest of the case of Conjecture \ref{conjmain}. Let $\overline{G}$ denote the complement of $G$. We {define $F_{a,b}$ as a star forest} of order $b+1$, which consists of a disjoint union of $\lfloor\frac{b+1}{a+1}\rfloor$ stars, ensuring that only one, if any, is non-isomorphic to $K_{1,a}$. Let $P^{\star}$ denote the Petersen graph. Let $uv$ be an edge in $G$ with the minimum $d(u) + d(v)$. The graph derived from $G$ by subdividing $uv$ a total of $k$ times is referred to as $S^{k}(G)$. Let $P_3=v_1v_2v_3$ and let $V_1\cup V_2\cup V_3$ be a partition of $V(K_{b-1})$, where $|V_i|=a_i$ for $i=1,2,3$. Then $F(a_1,a_2,a_3)$ denotes the graph obtained from $K_{b-1}$ and $P_3$ by adding  edges between each vertex of $V_i$ and $v_i$ for $i=1,2,3$.
\begin{thm}[\cite{Zhai-Lin}]\label{Zhai-Lin-1}
Given $b \geqslant a\geqslant 2$, $n-1=kb+t\ (0\leqslant t\leqslant b-1)$ and $\tau=\left\lfloor\frac{b+1}{a+1}\right\rfloor$,  let $G$ be the graph that attains the maximum $A_0$-spectral radius among all $K_{a,b}$-minor free graphs with $n$ vertices for large $n$.  
\begin{wst}
\item[{\rm (i)}]If $t\leqslant 2\tau-2$, then $G\cong K_{a-1} \vee ((k-1)K_b\cup F(a,0,b-1-a))$ for $t=\tau=2$ and $G\cong K_{a-1} \vee ((k-t)K_b\cup t\overline{F_{a,b}})$ otherwise;
\item[{\rm (ii)}]If $t\geqslant 2\tau-1$, then $G\cong K_{a-1} \vee ((k-1)K_b\cup \overline{P^{\star}})$ for $ t=2,\tau=1, b=8$ and $G\cong K_{a-1} \vee (kK_b\cup K_t)$ otherwise.
\end{wst}
\end{thm}

Notice that Zhai and Lin \cite{Zhai-Lin} also identified the $K_{1,b}$-minor free connected graphs having the maximum $A_0$-spectral radius. Up to now, the adjacency spectral extrema of $K_{a,b}$-minor free graphs are completely solved.


Besides focusing attention on extension of the confirmation for Conjecture \ref{conjmain}, our purpose here is to give a generalization of {Theorem~\ref{Zhai-Lin-1}} and some results of separated papers mentioned above. We give new and rather complete information about the relationship between the maximum $A_{\alpha}$-spectral radius and the structure of extremal $K_{a,b}$-minor free graphs. In details, we consider the extremal $K_{a,b}$-minor free graphs having the maximum $A_\alpha$-spectral radius for $1\leqslant a\leqslant b$ and  $0\leqslant \alpha<1$. Clearly, it suffices to consider the cases of $a=1,\ b\geqslant 3$ and $b\geqslant a \geqslant 2$. Our main results are described as follows.

\begin{thm}\label{main1}
Given $n\geqslant b \geqslant 3,$  let $G^{\star}$ be the connected graph that attains the maximum $A_\alpha$-spectral radius among all $K_{1,b}$-minor free graphs with $n$ vertices.
\begin{wst}
\item[{\rm (i)}]
If $n=b+1$, then $G^{\star}\cong \overline{F_{1,b}}$ for $\alpha\in[0,1)$.
\item[{\rm (ii)}]
If $n\neq b+1$, then
$G^{\star}\cong S^{n-b}(K_b)$ for $b=3, \ \alpha\in[0,1)$ or $b\geqslant 4, \alpha\in [\frac{2}{b+1},1)$.
\end{wst}
\end{thm}

\begin{thm}\label{main2}
For $2\leqslant a \leqslant b$, let $n-a+1=kb+t\ (0\leqslant t\leqslant b-1)$ and $\tau=\left\lfloor\frac{b+1}{a+1}\right\rfloor$, where $n$ is large enough. Denote by $G^*$ the graph that attains the largest $A_\alpha$-spectral radius among all $K_{a,b}$-minor free graphs with $n$ vertices for $\alpha\in[0,1)$.
\begin{wst}
\item[{\rm (i)}]If $t\leqslant 2(\tau-1)$, then $G^* \cong K_{a-1} \vee ((k-1)K_b\cup F(a,0,b-1-a))$ for $t=\tau=2$ and $G^* \cong K_{a-1} \vee ((k-t)K_b\cup t\overline{F_{a,b}})$ otherwise;
\item[{\rm (ii)}]If $t\geqslant 2\tau-1$, then $G^* \cong K_{a-1} \vee ((k-1)K_b\cup \overline{P^{\star}})$ for $ t=2,\tau=1, b=8$ and $G^* \cong K_{a-1} \vee (kK_b\cup K_t)$ otherwise.
\end{wst}
\end{thm}

Further on we need some notation. The set of neighbours (resp. non-neighbors) of a vertex $u$ in graph $G$ is denoted by $N_G(u)$\, (resp. $\overline{N}_G(u)$). Then $d_G(u)=|N_G(u)|$ is the degree of a vertex $u$ in a graph $G$. In particular, $\Delta(G)=\max\{d_G(u):\, u\in V(G)\}$. Let $G[V_1]$ denote the subgraph induced by $V_1\subseteq V(G)$. The \textit{distance} between two vertices $u$ and $v$, written as $d_G(u,v)$, is the length of a shortest path joining them. Let $E_G(A,B)=\{uv\in E(G):\, u\in A, v\in B\}$, and consequently, $e_G(A,B)=|E_G(A,B)|$. We denote the graph resulting from $G$ by removing vertex $v \in V(G)$ as $G-v$ and the graph resulting from $G$ after the deletion of edge $uv\in E(G)$ as $G-uv$. This notation is easily adapted for scenarios where multiple vertices or edges are removed. In a similar vein, $G + uv$ is derived from $G$ by the addition of edge $uv\not\in E(G)$. 

The remainder of this paper is organized as follows: In the next section, we recall some important known results and give some structure lemmas. In Section 3, we give the proof of Theorem \ref{main1}.
In Section 4, we give the proof of Theorem \ref{main2}. In the last section, some concluding remarks are given.
\section{\normalsize Preliminaries}
We give certain lemmas in this section that are useful for those that follow. Assume that $V=\{1, \ldots, n\}$ is the index of a real matrix $B$, with rows and columns. Assume that partition $V=V_{1} \cup \cdots \cup V_{k}$ enables $B$ to be expressed as 
$$
B = \begin{pmatrix}        B_{11}      &  \cdots  &  B_{1k}\\
                            \vdots     &  \ddots  & \vdots  \\
                            B_{k1}    & \cdots & B_{kk}
                    \end{pmatrix}
$$
where $B_{i j}$ represents the submatrix (block) of $B$ which is formed of the rows in $V_{i}$ and the columns in $V_{j}$. Denote by $c_{i j}$ the $B_{i j}$'s average row sum. In that case, matrix $C=(c_{i j})$ is referred to as $B$'s {\it quotient matrix}. For each $B_{i j}$, if the row sum is a constant, we call the partition {\it equitable} and $C$ the {\it equitable quotient matrix} of $B$.
\begin{lem}[\cite{0009}]\label{equit}
Given a nonnegative matrix $B,$  let $C$ represent the equitable quotient matrix of $B$. It holds true that the spectral radii of $C$ and $B$, denoted as  $\rho(C)$ and $\rho(B)$ respectively, are equal.
\end{lem}
For $\alpha\in[0,1)$, if graph $G$ is connected, then $A_{\alpha}(G)$ is an irreducible nonnegative matrix \cite{Niki1}. By {\rm Perron-Frobenius Theorem}, $A_{\alpha}(G)$ has a positive eigenvector ${\bf x}$ corresponding to $\lambda_\alpha(G)$. We call it the {\it Perron vector} of $A_{\alpha}(G)$. To simplify matters, we represent the coordinate of ${\bf x}$ with $x_v$, which aligns with vertex $v \in V(G)$. A straightforward calculation yields the following results
\begin{equation}\label{x}
\lambda_\alpha(G) x_v=\alpha d(v) x_v+(1-\alpha) \sum_{uv \in E(G)} x_u
\end{equation}
for all $v \in V(G)$.
\begin{lem}[\cite{Niki1}]\label{lambda}
Assume that $G$ contains no isolated vertices and $\alpha \in[0,1]$. Then
\begin{wst}
\item[{\rm (i)}]
$\lambda_\alpha(G) \leqslant\Delta(G)$, where the equality holds if and only if $G$ is $\Delta(G)$-regular;
\item[{\rm (ii)}]
$
\lambda_0(G)\leqslant \lambda_\alpha(G) \leqslant \max\limits_{v\in V(G)}\left\{\alpha d(v)+\frac{1-\alpha}{d(v)}\sum\limits_{uv\in E(G)}d(u)\right\}.
$
\end{wst}
\end{lem}

\begin{lem}[\cite{Niki1}]\label{sub}
Given a connected graph $G,$  if $H$ is a subgraph of $G$, then $\lambda_\alpha(H)\leqslant \lambda_\alpha(G)$, where $\alpha\in[0,1)$. In particular, if $H$ is a proper subgraph of $G$, then $\lambda_\alpha(H)<\lambda_\alpha(G)$.
\end{lem}
\begin{lem}[\cite{Niki2}]\label{trans1}
Given a connected graph $G,$ denote by ${\bf x}$ the Perron vector of $A_{\alpha}(G)$ for $\alpha\in[0,1)$. For $u, v \in V(G)$, assume that $v_1,\ldots,v_s \in N(v) \backslash \left(N(u)\cup\{u\}\right)$, where $1 \leqslant s \leqslant d(v)$. Let $G'=G-\{vv_i| 1 \leqslant i \leqslant s\}+\{uv_i: 1 \leqslant i \leqslant s\}$. If $x_{u}\geqslant x_{v}$, then $\lambda_\alpha(G')>\lambda_\alpha(G)$.
\end{lem}

In this paper, except for the Rayleigh quotient, double eigenvectors skill can also be used to compare the spectral radii of graphs, as described in the following lemma.
\begin{lem}[\cite{Lidan}]\label{xy}
Given two graphs $G$ and $H$ of equal order, if the Perron vectors corresponding to the matrices $A_{\alpha}(G)$ and $A_{\alpha}(H)$ are denoted as ${\bf x}$ and ${\bf y},$ respectively, we have
$$
{\bf x}^{T}(A_{\alpha}(G)){\bf y}=\sum_{uv\in E(G)}[\alpha (x_u y_u+x_v y_v)+(1-\alpha)(x_u y_v+x_v y_u)]
$$
and
$$
{\bf x}^{T} {\bf y}(\lambda_\alpha(H)-\lambda_\alpha(G))={\bf x}^{T}(A_{\alpha}(H)-A_{\alpha}(G)) {\bf y}.
$$
\end{lem}
The following lemma gives an upper bound of the edge number of $K_{1, b}$-minor free connected graph.
\begin{lem}[\cite{Chu,Ding}]\label{edge1b}
For $b \geqslant 3$ and $n \geqslant b+2$, let $G$ be a connected graph with $|V(G)|=n$ and $K_{1, b}$-minor free. Then we have $e(G) \leqslant {b \choose 2}+n-b$, and this holds as the optimal condition for $n, b$.
\end{lem}
Let $\tau=\left\lfloor\frac{b+1}{a+1}\right\rfloor$ and $\omega=\min\{a,\left\lfloor\frac{b+1}{2}\right\rfloor\}$. We say a graph $G$ has {\it the $(a,b)$-property}, if $G$ is $K_{r,s}$-minor free for any two positive integers $r,s$ with $r+s=b+1$ and $r \leqslant \omega$.
\begin{lem}[{\cite{Zhai-Lin}}]\label{b+21}
Given a connected graph $G$ with order $b+2$  having the $(a,b)$-property,  let $\tau\leqslant 2$. Then we have $e(G)\leqslant {b\choose 2}+2$, and we have $G$ is isomorphic to some $F(a_1,a_2,a_3)$ or $\overline{P^{\star}}$ if equality holds.
\end{lem}

\begin{lem}[\cite{Zhai-Lin}]\label{eb+1}
Given an edge-maximal connected graph $G$ of oder $b+1$ with the $(a,b)$-property,  then $e(G)={b\choose 2}+\tau-1$. Furthermore, $\overline{G}$ is a forest with $\tau$ components.
\end{lem}
\begin{lem}[\cite{Zhai-Lin}]\label{cb+1}
A connected graph $G$ with $b+1$ vertices has the $(a,b)$-property if and only if every component of $\overline{G}$ contains no fewer than $\omega +1$ vertices.
\end{lem}
\begin{lem}[{\cite{Zhai-Lin}}]\label{abgraph}
All the graphs $F(\omega, 0,b-1-\omega)$, $\overline{F_{a,b}}$ and $\overline{P^{\star}}$ have the $(a,b)$-property.
\end{lem}
A vertex $v\in V (G)$ is called a \textit{dominating vertex} of $G$ if $d_G(u)=|V(G)|-1$. A \textit{clique dominating set} is a set of dominating vertices. The next lemma gives {\b a} characterization of graphs with a clique dominating set.
\begin{lem}[\cite{Zhai-Lin}]\label{ab}
Let $G$ be a graph with a clique dominating set $S$ of size $a-1$. Then $G$ is $K_{a, b}$-minor free if and only if $G-S$ has the $(a,b)$-property.
\end{lem}
\begin{lem}\label{Mm}
Assume that graph $G$ has a clique dominating set $S$. Let ${\bf x}$ be the Perron vector of $A_{\alpha}(G)$.
Denote some components of $G-S$ by $\mathcal{F}$ and assume that $\Delta(\mathcal{F})<c$ (a constant).
Define
$$
X_s=\sum_{v\in S}x_v,\ X_M=\max_{v\in V(\mathcal{F})}x_v,\ X_m=\min_{v\in V(\mathcal{F})}x_v.
$$
\begin{wst}
\item[{\rm (i)}]One has
$
X_m\geqslant\frac{(1-\alpha)X_s}{\lambda_\alpha(G^*)-\alpha|S|},\
X_M<\frac{(1-\alpha)X_s}{\lambda_\alpha(G^*)-\alpha|S|-c}.
$
\item[{\rm (ii)}]
Given two constants $p,\ q$ with $p>q$, if $|V(G)|=n$ is large enough and $|S|$ is a constant, then $pX_m>qX_M$ and $pX_m^2>qX_M^2$. Furthermore, for two vertices $u,v\in V(\mathcal{F})$ with $d_{\mathcal{F}}(u)>d_{\mathcal{F}}(v)$, we have $x_u>x_v$.
\end{wst}
\end{lem}
\begin{proof}
In view of \eqref{x}, for each vertex $v\in V(\mathcal{F})$, we obtain
\[\label{F}
\lambda_\alpha(G)x_v=\alpha(|S|+d_{\mathcal{F}}(v))x_v+(1-\alpha)(X_s+\sum_{uv\in E(\mathcal{F})}x_u).
\]

{\rm (i)}\
By the definitions of $X_m$ and $X_M$, we have
$$
\alpha|S|X_m+(1-\alpha)X_s\leqslant \lambda_\alpha(G)X_m\leqslant
\lambda_\alpha(G)X_M<\alpha(|S|+c)X_M+(1-\alpha)(X_s+cX_M),
$$
and so
\[\label{X}
X_m\geqslant\frac{(1-\alpha)X_s}{\lambda_\alpha(G)-\alpha|S|},\
X_M<\frac{(1-\alpha)X_s}{\lambda_\alpha(G)-\alpha|S|-c}.
\]

{\rm (ii)}\
By Lemmas \ref{lambda}, \ref{sub} and $\Delta(G)=n-1$, we have $\lambda_\alpha(G)\geqslant\lambda_\alpha(K_{1,n-1})\geqslant\lambda_0(K_{1,n-1})={\sqrt{n-1}}$.
Note that $n$ is large enough and $p,q,c,|S|$ are all constants. By {\rm (i)},
\begin{align*}
pX_m-qX_M&>(1-\alpha)X_s\left(\frac{p}{\lambda_\alpha(G)-\alpha|S|}
-\frac{q}{\lambda_\alpha(G)-\alpha|S|-c}\right)\\
&=(1-\alpha)X_s\frac{(p-q)(\lambda_\alpha(G^*)-\alpha|S|)-cp}{(\lambda_\alpha(G)-\alpha|S|-c){(\lambda_\alpha(G)-\alpha|S|)}}\\
&>0.
\end{align*}
Similarly, we have
\begin{align*}
pX_m^2-qX_M^2&>(1-\alpha)^2X_s^2\left(\frac{p}{(\lambda_\alpha(G)-\alpha|S|)^2}
-\frac{q}{(\lambda_\alpha(G)-\alpha|S|-c)^2}\right)\\
&>(1-\alpha)^2X_s^2\left(\frac{\sqrt{p}}{(\lambda_\alpha(G)-\alpha|S|)^2}
-\frac{\sqrt{q}}{(\lambda_\alpha(G)-\alpha|S|-c)^2}\right)\\
&>0.
\end{align*}

For two vertices $u,v\in V(\mathcal{F})$, by \eqref{F} we have
\begin{align*}
\lambda_\alpha(G)x_u\geqslant \alpha(|S|+d_{\mathcal{F}}(u))X_m+(1-\alpha)(X_s+d_{\mathcal{F}}(u)X_m),\\
\lambda_\alpha(G)x_v\leqslant \alpha(|S|+d_{\mathcal{F}}(v))X_M+(1-\alpha)(X_s+d_{\mathcal{F}}(v)X_M).
\end{align*}
Note that $d_{\mathcal{F}}(u)>d_{\mathcal{F}}(v)$.
Then
\begin{align*}
\lambda_\alpha(G)(x_u-x_v)&\geqslant d_{\mathcal{F}}(u)X_m-d_{\mathcal{F}}(v)X_M>0,
\end{align*}
as desired.
\end{proof}
Let $X=\left(x_1, x_2, \ldots, x_n\right)^T$ and $Y=\left(y_1, y_2, \ldots, y_n\right)^T$ in $ \mathbb{R}^n$, with their components arranged in decreasing order. We say that $X$ is \textit{majorized} by $Y$, denoted by $X \prec Y$, if $\sum_{i=1}^k x_{i}$ is less than or equal to $\sum_{i=1}^k y_{i}$ for all $k=1,2, \ldots, n-1$, and $\sum_{i=1}^n x_i=\sum_{i=1}^n y_i$.
\begin{lem}[\cite{Niki2}]\label{major2}
Let $X, Y, Z\in \mathbb{R}^n$ with their components in decreasing order. If $X \prec Y$, we have $X^T \cdot Z \leqslant Y^T \cdot Z$.
\end{lem}

\section{\normalsize The Proof of Theorem \ref{main1}.}
In this section, we give the proof of  Theorem \ref{main1}, which determines the $n$ vertex $K_{1,b}$-minor free graph $G^{\star}$ having the maximum $A_\alpha$-spectral radius.

By Lemma \ref{lambda},  $\lambda_\alpha(G^{\star})\leqslant \Delta(G^{\star})\leqslant b-1$ for $\alpha\in[0,1)$.
Let $K_b^e$ represent the graph derived from $K_b$ by removing an edge $e=uv$ and then adding two pendant edges, say $uu',\ vv'$. The next lemma gives a lower bound of $\lambda_\alpha(G^{\star})$.
\begin{lem}\label{updown}
For $\alpha\in[0,1)$, we have $\lambda_\alpha(G^{\star})> b-1-\frac{2(1-\alpha)}{b-1}\geqslant b-2+\alpha$.
\end{lem}
\begin{proof}
By $b\geqslant 3$ and $\alpha\in[0,1)$, one has
$$b-1-\frac{2(1-\alpha)}{b-1}-(b-2+\alpha)=\frac{(1-\alpha)(b-3)}{b-1}\geqslant 0.$$
So, $b-1-\frac{2(1-\alpha)}{b-1}\geqslant b-2+\alpha$.
Notice that $S^{n-b}(K_b)$ is an $n$-vertex $K_{1,b}$-minor free graph. By the definition of $G^{\star}$, we get $\lambda_\alpha(G^{\star})\geqslant \lambda_\alpha(S^{n-b}(K_b))$.
Hence, we just need to show that
$\lambda_\alpha(S^{n-b}(K_b))>b-1-\frac{2(1-\alpha)}{b-1}$.

For $n=b+1$, the quotient matrix of $A_\alpha(S^1(K_b))$ is given as
$$
B_1=\begin{pmatrix}        2\alpha      & 2(1-\alpha)  & 0\\
                       (1-\alpha)   &  (b-1)\alpha &(b-2)(1-\alpha) \\
                       0            & 2(1-\alpha)  & (b-1)\alpha+(b-3)(1-\alpha)
                    \end{pmatrix}.
$$
By a direct calculation,  the characteristic polynomial of $B_1$ is
$$
f_1(x)=x^3-(\alpha b+b+3\alpha-3)x^2+(\alpha b^2+(2 \alpha^2+2 \alpha-2) b+2 \alpha^2-7 \alpha+2)x-2 \alpha^2 b^2+(2 \alpha^2+2) b-4 \alpha^2+8 \alpha-6.
$$
From Lemma \ref{equit}, $\lambda_\alpha(S^1(K_b))$ is the largest real root of $f_1(x)$.
On the other hand, we have
$$
f_1\left(b-1-\frac{2(1-\alpha)}{b-1}\right)=\frac{-4 (1-\alpha)^2 (b-2)((2-\alpha) b^2-4b+3\alpha )}{(b-1)^3}.
$$
By $\alpha\in[0,1)$, we get $\frac{2}{2-\alpha}<3\leqslant b$ and so
$(2-\alpha) b^2-4b+3\alpha \geqslant 6(1-\alpha)>0.$
Together with $b-2>b-1>0$ and $1-\alpha>0$.
We have $f_1\left(b-1-\frac{2(1-\alpha)}{b-1}\right)<0$.
Thus,
$\lambda_\alpha(S^1(K_b))>b-1-\frac{2(1-\alpha)}{b-1}$.

For $n\geqslant b+2$, we know that $K_b^e$ is a proper subgraph of $S^{n-b}(K_b)$.
By Lemma \ref{sub}, $\lambda_\alpha(S^{n-b}(K_b))>\lambda_\alpha(K_b^e)$.
Similar to the discussion for $n=b+1$, the largest real root of
$$
f_2(x)=x^3-(\alpha b+b+2\alpha-3)x^2+(\alpha b^2+(\alpha^2+\alpha-2)b+2 \alpha^2-6 \alpha+3)x-\alpha^2b^2+(\alpha^2+1)b-2 \alpha^2+4 \alpha-3
$$
is $\lambda_\alpha(K_b^e)$.
By some calculations, we have
$$
f_2\left(b-1-\frac{2(1-\alpha)}{b-1}\right)=\frac{-2 (1-\alpha)^2 (b-2) ((3-\alpha)b^2-6b+5\alpha-1)}{(b-1)^3}.
$$
Since $\alpha\in[0,1)$, we have $\frac{3}{3-\alpha}<3\leqslant b$ and so
$(3-\alpha)b^2-6b+5\alpha-1\geqslant 4(2-\alpha)>0.$
Note that $b-2>b-1>0$ and $1-\alpha>0$.
Then $f_2\left(b-1-\frac{2(1-\alpha)}{b-1}\right)<0$. So, $\lambda_\alpha(K_b^e)>b-1-\frac{2(1-\alpha)}{b-1}$,
i.e.,
$
\lambda_\alpha(S^{n-b}(K_b))>\lambda_\alpha(K_b^e)>b-1-\frac{2(1-\alpha)}{b-1}.
$
\end{proof}
Note that $G^{\star}$ is connected. Let ${\bf x}$ be the Perron vector of $G^{\star}$ and $x_{v^*}=\max\limits_{v\in V(G^{\star})}x_v$.
We get the next lemma.
\begin{lem}\label{U}
Let $U=N_{G^{\star}}(v^*)$. Then $|U|=b-1$.
\end{lem}
\begin{proof}
Since $G^{\star}$ is $K_{1,b}$-minor free, it is easy to get $|U|\leqslant \Delta(G^{\star})\leqslant b-1$.
Suppose that $|U|\leqslant b-2$.
By \eqref{x} and $x_{v^*}=\max\limits_{v\in V(G^{\star})}x_v$, we have
$$
\lambda_\alpha(G^{\star})x_{v^*}\leqslant \alpha|U|x_{v^*}+(1-\alpha)|U|x_{v^*}=|U|x_{v^*}\leqslant (b-2)x_{v^*}.
$$
Thus, $\lambda_\alpha(G^{\star})\leqslant|U|\leqslant b-2\leqslant b-2+\alpha$, which contradicts
$\lambda_\alpha(G^{\star})>b-2+\alpha$ (see Lemma \ref{updown}).
So, $|U|=b-1$.
\end{proof}
Given a partition of $U:=U_1\cup U_2$, where $U_1=\{v\in U|d_{U}(v)=b-2\}$.
Define $W=V(G^{\star})\backslash(U\cup \{v^*\})$ and $N^2(v^*)=\{w\in W|d_{G^{\star}}(v^*,w)=2\}$. We obtain the following lemma.
\begin{lem}\label{vertex}
For $v\in U\cup W$, we have $d_W(v)=0$ if $v\in U_1$, $d_W(v)\leqslant 1$ if $v\in U_2\cup N^2(v^*)$, and $d_W(v)\leqslant 2$ otherwise.
\end{lem}
\begin{proof}
For any vertex $v\in U_1$, it is easy to get $d_{G^{\star}}(v)=1+d_U(v)+d_W(v)=1+b-2+d_W(v)\leqslant \Delta(G^{\star})\leqslant b-1$. So, $d_W(v)=0$.

Suppose that there is a vertex $v\in U_2$ such that $d_W(v)\geqslant 2$.
In that case, $G^{\star}$ must have a double star as a subgraph with $b$ leaves and the central vertices are $v^*,\, v$. Since $G^{\star}$ is $K_{1,b}$-minor free, we get a contradiction.

For any vertex $v\in W$, as $G^{\star}$ is connected, there exists a shortest path, say $P$, connecting $v^*$ and $v$. Then $|N_{W\cap V(P)}(v)|\leqslant 1$. 
Clearly, if we contract the path $P$ to a new vertex $u$, then the resulting graph contains a star with central vertex $u$ and $b-2+|N_{W\backslash V(P)}(v)|$ leaves. Note that $G^{\star}$ is a $K_{1,b}$-minor free graph. Then $b-2+|N_{W\backslash V(P)}(v)|\leqslant b-1$ and so $|N_{W\backslash V(P)}(v)|\leqslant 1$. Then $d_W(v)=|N_{W\backslash V(P)}(v)|+|N_{W\cap V(P)}(v)|\leqslant 2$. In particular, if $v\in N^2(v^*)$, then $|N_{W\cap V(P)}(v)|=0$ and we have $d_W(v)\leqslant 1$.
\end{proof}
From Lemma \ref{vertex}, we know that $d_W(v)\leqslant 1$ for all vertices $v\in U_2$.  Hence, we obtain the next result.
\begin{cor}\label{11}
$|U_2|\geqslant |N^2(v^*)|$.
\end{cor}
\begin{lem}\label{trans1b}
Let $uv\in E(G^{\star})$ and $N_{G^{\star}}(u)\cap N_{G^{\star}}(v)=R$ with $d_{G^{\star}}(u)=d_{G^{\star}}(v)$.
Assume that $N_{G^{\star}}(u)\backslash (R\cup\{v\})=\{u_1\},\ N_{G^{\star}}(v)\backslash (R\cup\{u\})=\{v_1\}$.
If $x_{u}\geqslant x_v$, then $x_{u_1}\geqslant x_{v_1}$.
\end{lem}
\begin{proof}
By $\lambda_\alpha(G^{\star}){\bf x}=A_\alpha(G^{\star}){\bf x}$,  one sees
\begin{align*}
\lambda_\alpha(G^{\star})x_u&=\alpha d_{G^{\star}}(u)x_u+(1-\alpha)(x_v+x_{u_1}+\sum_{w\in R}x_w),\\
\lambda_\alpha(G^{\star})x_v&=\alpha d_{G^{\star}}(v)x_v+(1-\alpha)(x_u+x_{v_1}+\sum_{w\in R}x_w).
\end{align*}
So,
$\big(\lambda_\alpha(G^{\star})-\alpha d_{G^{\star}}(u)+(1-\alpha)\big)(x_u-x_v)=(1-\alpha)(x_{u_1}-x_{v_1})$.

Note that $d_{G^{\star}}(u)\leqslant \Delta(G^{\star})\leqslant b-1$ and $\lambda_\alpha(G^{\star})>b-2+\alpha$, we have
$$
\lambda_\alpha(G^{\star})-\alpha d_{G^{\star}}(u)+(1-\alpha)
\geqslant\lambda_\alpha(G^{\star})-(b-1)\alpha+(1-\alpha)
>b-2+\alpha-b\alpha+1
=(b-1)(1-\alpha)>0.
$$

Notice that $x_{u}\geqslant x_v$ and $1-\alpha>0$. Hence, $x_{u_1}-x_{v_1}\geqslant(b-1)(x_u-x_v)\geqslant 0$.
\end{proof}
Let $u^*\in V(G^{\star})$ such that
\begin{align*}
\alpha d_{G^{\star}}(u^*)+\frac{1-\alpha}{d_{G^{\star}}(u^*)} \sum_{v \in N_{G^{\star}}(u^*)} d_{G^{\star}}(v)
=\max\limits_{u \in V(G^{\star})} \left\{\alpha d_{G^{\star}}(u)+\frac{1-\alpha}{d_{G^{\star}}(u)} \sum_{v \in N_{G^{\star}}(u)} d_{G^{\star}}({v})\right\} .
\end{align*}
By Lemma \ref{lambda}, 
\[\label{ustar}
\lambda_\alpha(G^{\star}) \leqslant \alpha d_{G^{\star}}(u^*)+\frac{1-\alpha}{d_{G^{\star}}(u^*)} \sum_{v \in N_{G^{\star}}(u^*)} d_{G^{\star}}(v).
\]

In the next lemma, we give some information on the local structure of {$G^{\star}$}.
\begin{lem}\label{ustardegree}
Either $(b-1,b-1, \ldots, b-1, b-1)$ or $(b-2,b-1, \ldots, b-1, b-1)$ is the non-decreasing degree sequence of $N_{G^{\star}}(u^*)$. Moreover, $d_{G^{\star}}(u^*)=b-1$ if $\alpha\in[\frac{2}{b+1},1)$.
\end{lem}
\begin{proof}
Suppose that there exists a vertex $v\in N_{G^{\star}}(u^*)$ such that $d_{G^{\star}}(v)\leqslant b-3$. By \eqref{ustar} and Lemma~\ref{updown}, one has
\begin{align*}
b-1-\frac{2(1-\alpha)}{b-1}<\lambda_\alpha(G^{\star})
&\leqslant \alpha d_{G^{\star}}(u^*)+\frac{1-\alpha}{d_{G^{\star}}(u^*)} \sum_{v \in N_{G^{\star}}(u^*)} d_{G^{\star}}(v)\\
&\leqslant \alpha d_{G^{\star}}(u^*)+\frac{1-\alpha}{d_{G^{\star}}(u^*)}\big( (d_{G^{\star}}(u^*)-1)\Delta(G^{\star})+d_{G^{\star}}(v)\big)\\
&\leqslant \alpha d_{G^{\star}}(u^*)+\frac{1-\alpha}{d_{G^{\star}}(u^*)}\big( (b-2)(b-1)+(b-3)\big).
\end{align*}
By some calculations, we have
\begin{align}\label{u*}
\frac{(d_{G^{\star}}(u^*)-b+1)(2+((b-1)d_{G^{\star}}(u^*)-2)\alpha)}{d_{G^{\star}}(u^*)(b-1)}>0.
\end{align}
Notice that $d_{G^{\star}}(u^*)\geqslant 1,\ b\geqslant 3$ and $\alpha\in[0,1)$. Then $d_{G^{\star}}(u^*)>b-1=\Delta(G^{\star})$, a contradiction. Hence, $d_{G^{\star}}(v)\geqslant b-2$ for any vertex $v\in N_{G^{\star}}(u^*)$.

Suppose to the contrary that there exist at least two vertices $v,w\in N_{G^{\star}}(u^*)$ such that $d_{G^{\star}}(v)=d_{G^{\star}}(w)=b-2$. By Lemma~\ref{updown} and \eqref{ustar} , one has
\begin{align*}
b-1-\frac{2(1-\alpha)}{b-1}<\lambda_\alpha(G^{\star})
&\leqslant \alpha d_{G^{\star}}(u^*)+\frac{1-\alpha}{d_{G^{\star}}(u^*)} \sum_{v \in N_{G^{\star}}(u^*)} d_{G^{\star}}(v)\\
&\leqslant \alpha d_{G^{\star}}(u^*)+\frac{1-\alpha}{d_{G^{\star}}(u^*)}\big( (d_{G^{\star}}(u^*)-2)\Delta(G^{\star})+d_{G^{\star}}(v)+d_{G^{\star}}(w)\big)\\
&\leqslant \alpha d_{G^{\star}}(u^*)+\frac{1-\alpha}{d_{G^{\star}}(u^*)}\big( (b-3)(b-1)+2(b-2)\big).
\end{align*}
By some calculations, we may deduce \eqref{u*}. Similar to the former discussion, we get a contradiction.

To complete the proof, we need to show that $d_{G^{\star}}(u^*)=b-1$ if $\alpha\in[\frac{2}{b+1},1)$. Recall $d_{G^{\star}}(u^*)\leqslant \Delta(G^{\star})=b-1$.
By Lemma \ref{updown} and \eqref{ustar}, we have
\begin{align*}
b-1-\frac{2(1-\alpha)}{b-1}< \lambda_\alpha(G^{\star})
&\leqslant \alpha d_{G^{\star}}(u^*)+\frac{1-\alpha}{d_{G^{\star}}(u^*)} \sum_{v \in N_{G^{\star}}(u^*)} d_{G^{\star}}(v) \\
& \leqslant \alpha d_{G^{\star}}(u^*)+(1-\alpha)(b-1).
\end{align*}
This implies that $d_{G^{\star}}(u^*)>b-2+\frac{\alpha(b+1)-2}{\alpha(b-1)}$. If $\alpha\in [\frac{2}{b+1},1)$, then $d_{G^{\star}}(u^*)=b-1$.
\end{proof}

\subsection{\normalsize The proof of Theorem~\ref{main1}(i)}
In this subsection, we present the proof of Theorem~\ref{main1}(i). If $n=b+1$, then $|W|=|V(G^{\star})\backslash(U\cup \{v^*\})|=1$. Let $W=\{w\}$. We give the following lemma.
\begin{lem}\label{U2}
    $d_W(v)=1$ and $d_U(v)=b-3$ for any $v\in U_2$.
\end{lem}
\begin{proof}
By Lemma~\ref{vertex}, we only need to show that $d_W(v)\not=0$. Otherwise, we have $vw\notin E(G^{\star})$.
Let $G'=G^{\star}+vw$.
Since $d_{G^{\star}}(v),\ d_{G^{\star}}(w)\leqslant b-2$, we get $\Delta(G')=\Delta(G^{\star})=b-1$.
Notice that that $n=b+1$. We have $G'$ is a $K_{1,b}$-minor free graph.
By $G^{\star}$ is a proper subgraph of $G'$ and Lemma \ref{sub}, we have $\lambda_\alpha(G')>\lambda_\alpha(G^{\star})$, a contradiction.
Thus, $d_W(v)=1$ for each $v\in U_2$.

Suppose that there exists a vertex $v\in U_2$ such that $d_U(v)\neq b-3$.
By the definition of $U_2$, we have $d_U(v)\leqslant b-4$.
Choose a vertex $u\in U$ such that $uv\notin E(G^{\star})$.
Then $d_U(u)\leqslant b-3$ and so $u\in U_2$.
We claim that $d_U(u)=b-3$.
Otherwise, let $G'=G^{\star}+uv$.
Then $\Delta(G')=\Delta(G^{\star})=b-1$.
Notice that $n=b+1$. Then $G'$ is $K_{1,b}$-minor free.
By Lemma \ref{sub}, we have $\lambda_\alpha(G')>\lambda_\alpha(G^{\star})$, a contradiction.
Thus, $d_U(u)=b-3$ and so $d_{G^{\star}}(u)=1+d_U(u)+d_W(u)=1+b-3+1=b-1$.
So, $\overline{N}_{G^{\star}}(u)=\{v\}$.
Let $U\backslash\{u,v\}=R$.
Then $v^*u\in E(G^{\star})$,
$N_{G^{\star}}(v^*)\cap N_{G^{\star}}(u)=R$, $d_{G^{\star}}(v^*)=d_{G^{\star}}(u)$,
$N_{G^{\star}}(v^*)\backslash (R\cup\{u\})=\{v\}$, $N_{G^{\star}}(u)\backslash (R\cup\{v^*\})=\{w\}$.
By Lemma \ref{trans1b} and $x_{v^*}\geqslant x_u$, we have $x_v\geqslant x_w$.
Let $G'=G^{\star}-wu+vu$. Then $\Delta(G')=\Delta(G^{\star})=b-1$.
Since $n=b+1$, $G'$ is $K_{1,b}$-minor free.
Together with Lemma~\ref{trans1} and $x_v\geqslant x_w$, we get $\lambda_\alpha(G')>\lambda_\alpha(G^{\star})$, which contradicts the maximality of $\lambda_\alpha(G^{\star})$.
\end{proof}

\begin{proof}[\bf The proof of Theorem \ref{main1} {\rm (i)}]
According to Lemma \ref{U2}, for each $v\in U_2$, there is a unique non-neighbour of $v$ in $U_2$. Then $G^{\star}[U_2]\cong K_{|U_2|}-M,$ where $M$ is a perfect matching of $K_{|U_2|}$.
Thus, $|U_2|$ is even. 
We now consider the next two cases.

{\bf Case 1.}\ $n=b+1$ is even.
For a connected graph $G$, $\lambda_\alpha(G)\leqslant \Delta(G)$ where the equality holds if and only if $G$ is a $\Delta(G)$-regular graph (see Lemma \ref{lambda}).
Recall $d_{G^{\star}}(v^*)=|U|=b-1$.
For any vertex $v\in U$, by the definition of $U_1$ and Lemma \ref{U2}, we have $d_{G^{\star}}(v)=b-1$.
Thus, $G^{\star}$ is $\Delta(G)$-regular graph is equivalent to $d_{G^{\star}}(w)=b-1$.
Recall $N_{G^{\star}}(w)=N_{U}(W)=U_2$.
We have $|U_2|=b-1=|U|$, i.e., $U=U_2$.
Therefore, $\overline{G^{\star}}\cong M\cup\{v^*w\}=\frac{n}{2}K_2\cong F_{1,b}$.

{\bf Case 2.}\ $n=b+1$ is odd. Clearly, $b$ is even. Recall that $|U|=b-1$ and $|U_2|$ is even. Hence, $2\leqslant|U_2|\leqslant b-2$. By symmetry, we may set $x_1=x_v$ for all vertices $v\in U_1$ and $x_2=x_u$ for all vertices $u\in U_2$. Furthermore, we have $x_{v^*}=x_1$. Let $\lambda_\alpha^{\star}:=\lambda_\alpha(G^{\star})$. Then
\begin{align*}
\lambda_\alpha^{\star}x_{v^*}&=\alpha(b-1)x_{v^*}+(1-\alpha)(|U_1|x_{v^*}+|U_2|x_2),\\[3pt]
\lambda_\alpha^{\star}x_2&=\alpha(b-1)x_2+(1-\alpha)((|U_1|+1)x_{v^*}+(|U_2|-2)x_2+x_w),\\[3pt]
\lambda_\alpha^{\star}x_w&=\alpha|U_2|x_w+(1-\alpha)|U_2|x_2.
\end{align*}
Note that $|U|=|U_1|+|U_2|=b-1$. Let
\begin{align*}
h(x)
&={x}^3
-(\alpha|U_2|+b\alpha +\alpha+b-3){x}^2
+((b\alpha^2 +\alpha^2+b\alpha-3\alpha)|U_2|+b^2\alpha-\alpha -2b+2)x,\\
g(y)&=(1-\alpha)^2y^2+(b\alpha^2-1)(b-1)y,
\end{align*}
where $x\geqslant b-2+\alpha$ and $2\leqslant y\leqslant b-2$.
It is straightforward to check that $h(\lambda_\alpha^{\star})=g(|U_2|).$

By some calculations, we obtain
\begin{align*}
&h'(x)=
3{x}^2
-2(\alpha|U_2|+b\alpha +\alpha+b-3){x}
+(b\alpha^2 +\alpha^2+b\alpha-3\alpha)|U_2|+b^2\alpha-\alpha -2b+2,\\
&h''(x)=6{x}
-2(\alpha|U_2|+b\alpha +\alpha+b-3),\ h'''(x)=6>0.
\end{align*}
Then for $x\geqslant b-2+\alpha$, $h''(x)$ is a monotonically increasing function, (see Lemma \ref{updown}).
By $2\leqslant|U_2|\leqslant b-2$ and $b\geqslant 4$ ($b+1$ is odd),
we have
$$
h''(x)\geqslant h''(b-2+\alpha)=2((2-b-|U_2|)\alpha+2b-3)>2(b-1-|U_2|)>0.
$$
This implies that $h'(x)$ is a monotonically increasing function in $x$ for $x\geqslant b-2+\alpha$.
We have
\begin{align*}
h'(x)>h'(b-2+\alpha)
=&(1-\alpha)(b^2-(4-2\alpha)b-\alpha+2-\alpha(b-1)|U_2|)\\
\geqslant&(1-\alpha)((1 - \alpha) b^2 + (5\alpha-4)b- 3\alpha+2)\\
\geqslant&(1-\alpha)(2+\alpha)>0.
\end{align*}
Then $h(x)$ is a monotonically increasing function in $x>b-2+\alpha$.

By $\lambda_\alpha^{\star}>b-2+\alpha$, the maximality of $\lambda_\alpha^{\star}$ and $h(\lambda_\alpha^{\star})=g(|U_2|)$, we get $g(|U_2|)$ attains its maximum. By simple calculations, we have $g'(y)=2(1-\alpha)^2y+(b\alpha^2-1)(b-1).$ By $2(1-\alpha)^2>0$, then $\max\limits_{2\leqslant y\leqslant b-2}g(y)=\max\{g(2),\ g(b-2)\}.$ Since $g(b-2)-g(2)=(b-4)(\alpha b-1)^2\geqslant0$, we get $\max\limits_{2\leqslant y\leqslant b-2} g(y)=g(b-2).$ So, $|U_2|=b-2$ and $|U_1|=1$. Let $U_1=\{u\}$. Then $\overline{G^{\star}}=M\cup\{v^*w, uw\}=\frac{n-3}{2}K_2\cup P_3\cong F_{1,b}$.
\end{proof}

\subsection{\normalsize The proof of Theorem \ref{main1}{\rm (ii)}}
In this subsection, we give the proof of Theorem \ref{main1}(ii). We need the following two lemmas first.
\begin{lem}\label{n2v}
For $n\geqslant b+2$ and $\alpha\in[\frac{2}{b+1},1)$, we have $|N^2(v^*)|\geqslant 2$.
\end{lem}
\begin{proof}
By $|U|=b-1$ and $G^{\star}$ is connected, we have $|W|=n-b\geqslant 2$ and so $N^2(v^*)\neq\emptyset$.
Let $W=\{w_1,w_2,\ldots,w_{n-b}\}$.

Suppose that $N^2(v^*)=\{w_1\}$.
Let $v_1\in N_U(w_1)$. Then $v_1\in U_2$.
By Lemma \ref{vertex}, we have $d_W(w_1)=1$ and $G^{\star}[W]$ is a pendant path. Then we may set $N_W(w_1)=\{w_2\}$ and $G^{\star}[W]\cong w_1w_2\cdots w_{n-b}$.

If $v_1v\in E(G^{\star})$ for all vertices $v\in \overline{N}_{U_2}(w_1)$, then there is a double star with central vertices $v_1,\ w_1$ and $b$ leaves in $G^{\star}$, i.e., $G^{\star}$ contains a $K_{1,b}$-minor, a contradiction. Hence, there exists a vertex $v_2\in \overline{N}_{U_2}(w_1)$ such that $v_1v_2\notin E(G^{\star})$. We have $d_W(v_2)=0$, $d_U(v_2)\leqslant b-3$. Then $v_2\in U_2$ and $d_{G^{\star}}(v_2)\leqslant b-2$.

In order to complete our proof, we should characterize some local structure on $G^{\star}.$ We show the following claims.
\begin{claim}
$d_{G^{\star}}(w_1)\geqslant 3$.
\end{claim}
\begin{proof}
Suppose that $d_{G^{\star}}(w_1)=2$. We have $N_{G^{\star}}(w_1)=\{v_1,w_2\}$. Let $G'=G^{\star}+v_2w_{n-b}$. Then $G'$ is connected and $S^{n-b}(K_b)$ contains $G'$ as a subgraph. So, $G'$ is $K_{1,b}$-minor free.
By Lemma \ref{sub}, $\lambda_\alpha(G')>\lambda_\alpha(G^{\star})$, a contradiction.
\end{proof}

\begin{claim}
$b\geqslant 5$.
\end{claim}
\begin{proof}
Suppose that $b=4$, we have $|U|=b-1=3$. Let $U\backslash\{v_1,v_2\}=\{v_3\}$. Notice that $v_2w_1\notin E(G^{\star}), d_W(w_1)=1$ and $d_{G^{\star}}(w_1)\geqslant 3$. Consequently, $N_U(w_1)=\{v_1,v_3\}$. If $v_2v_3\in E(G^{\star})$, then $G^{\star}$ contains a double star with central vertices $v_3,w_1$ and $b$ leaves, i.e., $G^{\star}$ contains a $K_{1,b}$-minor, a contradiction. Thus, $v_2v_3\notin E(G^{\star})$.

Combing with $v_2v_1,\ v_2w_1\notin E(G^{\star})$ gives us $N_{G^{\star}}(v_2)=\{v^*\}$. From $v^*w_1\notin E(G^{\star})$, we have $G^{\star}[\{v^*, w_1, v_1,$ $ v_3\}]\subseteq K_4-v^*w_1$.
Let $G'=G^{\star}+v_2w_{n-b}$. Then $G'\subseteq S^{n-b}(K_b)$ and so $G'$ is still $K_{1,b}$-minor free.
By Lemma \ref{sub}, we get $\lambda_\alpha(G')>\lambda_\alpha(G^{\star})$, a contradiction.
\end{proof}

\begin{claim}\label{uin}
$u^*\notin W$.
\end{claim}
\begin{proof}
If not, by Lemma \ref{vertex}, we get $d_{G^{\star}}(w_i)\leqslant 2$ for $i=2,\ldots, n-b$. Together with Lemma \ref{ustardegree} and $b\geqslant 5$, we have $d_{G^{\star}}(u^*)=b-1\geqslant 4$.
Therefore, $u^*=w_1$. Recall $u^*w_2\in E(G^{\star})$. By Lemma \ref{ustardegree}, one has $d_{G^{\star}}(w_2)\geqslant b-2\geqslant 3$, a contradiction to $d_{G^{\star}}(w_2)\leqslant 2$. Thus, $u^*\notin W$.
\end{proof}

We come back to show Lemma~\ref{n2v}. In fact, we proceed by considering the following cases (based on Claim~\ref{uin}).

{\bf Case 1.}\ $u^*\notin U_2$. Then $u^*\in \{v^*\}\cup U_1$. We assume that $v^*=u^*$. By Lemma \ref{ustardegree}, $d_{G^{\star}}(v)\geqslant b-2$ for all vertices $v\in U$. Since {$v_2\in U_2$,} $d_{G^{\star}}(v_2)\leqslant b-2$. We obtain $d_{G^{\star}}(v_2)=b-2$. Hence, $\overline{N}_U(v_2)=\{v_1\}=U_2\setminus \{v_1,v_2\}$.

By $d_U(w_1)=d_{G^{\star}}(w_1)-1\geqslant 2$, we may choose $v_3\in N_{U}(w_1)$ such that $v_3\neq v_1$. If $v_3v\in E(G^{\star})$ for all $v\in \overline{N}_{U}(w_1)$, then $G^{\star}$ contains a double star with central vertices $v_3,w_1$ and $b$ leaves, a contradiction. Then we may set $v_4\in \overline{N}_{U}(w_1)$ and $v_4v_3\notin E(G^{\star})$. We have $d_W(v_4)=0$, $d_U(v_4)\leqslant b-3$ and so $d_{G^{\star}}(v_4)\leqslant b-2$. Notice that $\overline{N}_U(v_2)=\{v_1\}$. Then $v_2\neq v_4$. Recall $d_{G^{\star}}(v_2)=b-2,\ d_{G^{\star}}(v_4)\leqslant b-2$ and {$v_2,\ v_4\in N_{G^{\star}}(v^*)$,} we get a contradiction to Lemma \ref{ustardegree}.

{\bf Case 2.}\ $u^*\in U_2$. Note that $d_{G^{\star}}(u^*)=b-1$. Hence $w_1u^*\in E(G^{\star})$. We may assume that $v_1=u^*$. Then $d_{G^{\star}}(v_1)=b-1$ and so $d_U(v_1)=b-3$. As $v_1v_2\notin E(G^{\star})$, one has $\overline{N}_{U_2}(v_1)=\{v_2\}$. In view of the definition of $U_2$, $d_U(v)\leqslant b-3$ for any $v\in U_2\backslash \{v_1,v_2\}$. Notice that $w_1v_1\in E(G^{\star})$.
Hence, together with Lemma \ref{ustardegree}, we have $d_{G^{\star}}(w_1)\geqslant b-2$.

{\bf Case 2.1.}\ $d_{G^{\star}}(w_1)=b-2$. By Lemma \ref{ustardegree}, $d_{G^{\star}}(v)=b-1$ for $v\in U_2\backslash \{v_1,v_2\}$, and so $d_W(v)=1$, i.e., $vw_1\in E(G^{\star})$. Therefore, $N_{U_2}(w_1)=U_2\backslash \{v_2\}$.

Suppose that there exists a vertex $v'\in N_{U_2}(w_1)\backslash \{v_1\}$ such that $v'v_2\in E(G^{\star})$. Then $\overline{N}_{U_2}(v')\subseteq U_2\backslash \{v_2\}=N_{U_2}(w_1)$. Consequently, $G^{\star}$ contains a double star with central vertices $v',w_1$ and $b$ leaves, a contradiction. Therefore, $vv_2\notin E(G^{\star})$ for all $v\in N_{U_2}(w_1)\backslash\{v_1\}$. Notice that $d_{U}(v)=b-3$. Hence, $\overline{N}_{U_2}(v)=\{v_2\}$ for all  $v\in N_{U_2}(w_1)\backslash \{v_1\}$.

Let $R=U\backslash\{v_1,v_2\}$. We have $v^*v_1\in E(G^{\star})$, $N_{G^{\star}}(v^*)\cap N_{G^{\star}}(v_1)=R$, $d_{G^{\star}}(v^*)=d_{G^{\star}}(v_1)$, $N_{G^{\star}}(v^*) \backslash (R\cup\{v_1\})=\{v_2\}$, $N_{G^{\star}}(v_1) \backslash (R\cup\{v^*\})=\{w_1\}$. By Lemma \ref{trans1b} and $x_{v^*}\geqslant x_{v_1}$, one has $x_{v_2}\geqslant x_{w_1}$. Let $G'=G^{\star}-\sum_{v\in N_{U_2}(w_1)\backslash \{v_1\}}vw_1+\sum_{v\in N_{U_2}(w_1)\backslash \{v_1\}}vv_2.$ Then $G'[U\cup\{v^*\}]\cong K_b-v_1v_2$. Bearing in mind that $G'[W]$ is a pendant path $w_1\ldots w_{n-b}$ and $E_{G'}(U,W)=\{v_1w_1\}$, we obtain that $S^{n-b}(K_b)$ contains $G'$ as a proper subgraph.
So, $G'$ is $K_{1,b}$-minor free. By $x_{v_2}\geqslant x_{w_1}$ and Lemma \ref{trans1}, $\lambda_\alpha(G')>\lambda_\alpha(G^{\star})$, we get a contradiction.

{\bf Case 2.2.}\ $d_{G^{\star}}(w_1)=b-1$. As $N_W(w_1)=\{w_2\}$, we have $d_{U_2}(w_1)=b-2$. Note that $v_2\in U_2$ and $v_2w_1\notin E(G^{\star})$. Hence, $|\overline{N}_{U_2}(w_1)|\geqslant 1$ and $b-1=|U|\geqslant |U_2|=d_{U_2}(w_1)+|\overline{N}_{U_2}(w_1)|\geqslant b-1$. Thus, $U_2=U$ and $\overline{N}_U(w_1)=\{v_2\}$.

Suppose that we can find $v'\in N_U(w_1)$ satisfying that $v'v_2\in E(G^{\star})$. Together with $\overline{N}_U(w_1)=\{v_2\}$, one sees that $G^{\star}$ contains a double star with central vertices $v',w_1$ and $b$ leaves, a contradiction. So we have $vv_2\notin E(G^{\star})$ for all $v\in N_U(w_1)$.

Combing with $N_U(w_1)\cup \{v_2\}=U$ and $v_2w_1\notin E(G^{\star})$ gives us $N_{G^{\star}}(v_2)=\{v^*\}$. By Lemma \ref{ustardegree} and $\overline{N}_U(u^*)=\overline{N}_{U_2}(u^*)=\{v_2\}$,
one has $d_{G^{\star}}(u)\geqslant b-2$ for all $u\in U\backslash\{v_2\}$. Note that $N_U(w_1)=U\backslash\{v_2\}$. Consequently, $d_{G^{\star}}(u)\geqslant b-2$ for all $u\in N_U(w_1)$.

Suppose that there exists a vertex $v'\in N_{U}(w_1)$ such that $d_{G^{\star}}(v')=b-2$. Note that $v^*v',\ v'w_1\in E(G^{\star})$. We obtain $d_U(v')=b-4$ and so $|\overline{N}_U(v')|=2$. Recall $v_2\in \overline{N}_U(v')$. Then there exists another vertex $v''$ such that $v''\in (U\backslash \{v_2,v'\})\cap \overline{N}_U(v')$. Together with $v''v_2\notin E(G^{\star})$, one may deduce that $d_U(v'')\leqslant b-4$. Then $d_{G^{\star}}(v'')\leqslant b-2$. As $v',\ v''\in U\setminus \{v_2\}$, one has $v',\ v'' \in N_{G^{\star}}(u^*)$. By Lemma~\ref{ustardegree}, at most one vertex in $\{v',\, v''\}$ is of degree $b-2$, a contradiction. Thus, $d_{G^{\star}}(v)=b-1$ for $v\in N_{U}(w_1)$. Consequently, $G^{\star}[(U\backslash\{v_2\})\cup\{v^*,w_1\}]\cong K_b-v^*w_1$.

Let $G'=G^{\star}+v_2w_{n-b}$. Then $N_{G'}(v_2)=\{v^*,w_{n-b}\}$ and $G'[\{w_2,\ldots,w_{n-b},v_2\}]$ is a path. So, $G'\cong S^{n-b}(K_b)$, i.e., $G'$ is $K_{1,b}$-minor free. By Lemma \ref{sub}, we have $\lambda_\alpha(G')>\lambda_\alpha(G^{\star})$, a contradiction.
\end{proof}
\begin{lem}\label{u2w}
If $n\geqslant b+2$ and $\alpha\in[\frac{2}{b+1},1)$, then $d_W(v)=1$ for every $v\in U_2$.
\end{lem}
\begin{proof}
Suppose that there exists a vertex $v_1 \in U_2$ such that $d_W(v_1)\neq 1$. By Lemma~\ref{vertex} and the definition of $U_2$, we have $d_U(v_1)\leqslant b-3$ and $d_W(v_1) \leqslant 1$.
So, $d_W(v_1)=0$ and $d_{G^{\star}}(v_1)=d_U(v_1)+1\leqslant b-2$.

Combining with $n\geqslant b+2$ and Lemma \ref{n2v}, we have $|N^2(v^*)| \geqslant 2$. Suppose $|N_U(W)|=1$ and let $N_U(W)=\{u\}$. Then $d_W(u)=|N^2(v^*)|\geqslant 2$. Note that $u\in U_2$. Together with Lemma \ref{vertex}, we have $d_W(u) \leqslant 1$, a contradiction. Hence, $|N_U(W)|\geqslant 2$. By Lemma~\ref{vertex}, for all $w \in W$, we have $d_{G^{\star}}(w) =d_W(w) \leqslant 2$ if $w \in W\backslash N^2(v^*)$ and $d_W(w) \leqslant 1$ otherwise. Note that $|N^2(v^*)| \geqslant 2$ and $v_1w\notin E(G^{\star})$ for $w \in N^2(v^*)$. Hence, $d_U(w) \leqslant b-3$ and so $d_{G^{\star}}(w)=d_U(w)+d_W(w) \leqslant b-2$. Then $d_{G^{\star}}(w) \leqslant b-2$ for every $w \in W$.

By Lemma~\ref{ustardegree}, $d_{G^{\star}}(u^*)=b-1$. Therefore, $u^* \notin W \cup\{v_1\}$, and so $u^* \in(U_1 \cup\{v^*\}) \cup(U_2 \backslash\{v_1\})$. We proceed by considering the following two possible cases.

{\bf Case 1.}\ $u^* \in U_1 \cup\{v^*\}$. Clearly, $u^*v_1 \in E(G^{\star})$. Consequently, together with Lemma \ref{ustardegree}, one has $d_{G^{\star}}(v_1) \geqslant b-2$. Recall that $d_{G^{\star}}(v_1) \leqslant b-2$ and $d_W(v_1)=0$. Hence, $d_{G^{\star}}(v_1)=b-2$ and $d_U(v_1)=b-3$. Notice that $|U|=b-1$. So we have $|\overline{N_{U_2}}(v_1)|=1$. Then let $\overline{N_{U_2}}(v_1)=\{v_2\}$. Consequently, $v_2\in U_2$ and $v_2 \in N_{G^{\star}}(u^*)\backslash \{v_1\}$. By Lemma \ref{ustardegree}, we have $d_{G^{\star}}(v_2)=b-1$. In view of Lemma \ref{vertex} and by the definition of $U_2$, one has $d_U(v_2)\leqslant b-3$ and $d_W(v_2)\leqslant 1$. Therefore, $d_W(v_2)=1, d_U(v_2)=b-3,$ and so $N_U(v_2)=U \backslash \{v_1, v_2\}$. Let $N_W(v_2)=\{w_1\}$. Recall that $|N^2(v^*)| \geqslant 2$. Then there exists a vertex, say $w_2,$ in $N^2(v^*) \backslash\{w_1\}$. So we may assume  $v_3 \in N_U(w_2)$. Clearly, $v_3v_1, v_3v_2 \in E(G^{\star})$. Therefore, there is a double star with central vertices $v_2,v_3$ and $b$ leaves in $G^{\star}$, a contradiction.

{\bf Case 2.}\  $u^*\in U_2 \backslash\{v_1\}$. By Lemma \ref{vertex} and the definition of $U_2$, we have $d_U(u^*)\leqslant b-3$ and $d_W(u^*)\leqslant 1$. Notice that $d_{G^{\star}}(u^*)=b-1$. Hence, $d_W(u^*)=1$ and $d_U(u^*)=b-3$. So let $N_W(u^*)=\{w_1\}$. Consequently, $d_{G^{\star}}(w_1) \geqslant b-2$ (based on Lemma \ref{ustardegree}). Note that $w_1 \in W$. Hence, $d_{G^{\star}}(w_1) \leqslant b-2$, and so $d_{G^{\star}}(w_1)=b-2$.
By Lemma \ref{vertex}, one sees that $d_W(w_1)\leqslant 1$. Hence, $d_U(w_1)= d_{G^{\star}}(w_1)-d_W(w_1)\geqslant b-3$. By $|U \backslash\{v_1\}|=b-2$ and $d_W(v)\leqslant 1$ for all $v\in U_2$, we obtain
$|N^2(v^*)|\leqslant 2$. Together with $|N^2(v^*)|\geqslant 2$, one has $|N^2(v^*)|=2$.

Assume that $N^2(v^*)=\{w_1,\ w_2\}$. Then $b-1=|U|\geqslant |U_2|=|N_U(w_1)\cup N_U(w_2)\cup \{v_1\}|\geqslant b-1$. So, $U_2=U$, $d_{U_2}(w_1)=b-3$ and $d_U(w_2)=1$. Assume $N_U(w_2)=\{v_2\}$. Then we have $\overline{N}_U(w_1)=\{v_2,v_1\}$. Notice that $d_{G^{\star}}(w_1)=b-2$ and $w_1 u^* \in E(G^{\star})$. Hence, by Lemma \ref{ustardegree}, $d_{G^{\star}}(v)=b-1$ for all $v\in N_{G^{\star}}(u^*)\backslash \{w_1\}$.
Since $d_{G^{\star}}(v_1) \leqslant b-2$, one sees $v_1 \notin N_{G^{\star}}(u^*)$. Recall that $d_U(u^*)=b-3$. So we have $N_U(u^*)=U \backslash \{v_1, u^*\}$.

Suppose that $v_2v_1 \in E(G^{\star})$. Together with $v_2u^* \in E(G^{\star})$, we can find that there is a double star with central vertices $v_2,u^*$ and $b$ leaves in $G^{\star}$, a contradiction. Hence, $v_2v_1 \notin E(G^{\star})$. As $v_2\in N_{G^{\star}}(u^*)\backslash \{w_1\}$ and $d_U(w_2)=1$, we have $d_{G^{\star}}(v_2)=b-1$ and $d_U(v_2)=b-3$. Thus, $\overline{N}_U(v_2)=\{v_1\}$.

Together with Lemma \ref{ustardegree}, $N_U(u^*) \backslash\{v_2\}\subseteq N_{G^{\star}}(u^*)\backslash \{w_1\}$ and $d_{G^{\star}}(w_1)=b-2$, we may deduce $d_{G^{\star}}(v)=b-1$ for all $v \in N_U(u^*) \backslash\{v_2\}$.
Note that $N_{G^{\star}}(W)=U\backslash\{v_1\}$. Hence, $d_U(v)=d_{G^{\star}}(v)-2=b-3$, and so $|\overline{N}_{U_2}(v)|=1$. Suppose that there exists a vertex $v'\in N_U(u^*) \backslash\{v_2\}$ such that $v'v_1\in E(G^{\star})$.
Combining with $\overline{N}_U(w_1)=\{v_2,v_1\}$ and $d_W(w_1)=d_{G^{\star}}(w_1)-d_{U_2}(w_1)=1$ gives us that $G^{\star}$ contains a double star with central vertices $v',w_1$ and $b$ leaves, a contradiction.
Hence, for all $v \in N_U(u^*) \backslash\{v_2\}$, we have $vv_1 \notin E(G^{\star})$, i.e., $\overline{N}_{U_2}(v)=\{v_1\}$.

Let $U\backslash\{v_1, u^*\}=R$. Then it is easy to check that $v^*u^* \in E(G^{\star}), N_{G^{\star}}(v^*) \cap N_{G^{\star}}(u^*)=R$, $d_{G^{\star}}(v^*)=d_{G^{\star}}(u^*),$ $ N_{G^{\star}}(u^*) \backslash (R \cup\{v^*\})=\{w_1\}$ and $N_{G^{\star}}(v^*) \backslash(R \cup\{u^*\})=\{v_1\}$. By Lemma \ref{trans1b} and $x_{v^*} \geqslant x_{u^*}$, we obtain $x_{v_1} \geqslant x_{w_1}$.

By Lemma \ref{vertex} and $G^{\star}$ is connected, $G^{\star}[W]$ is a union of two paths at most. Notice that $N_{U}(w_1)=U\backslash\{v_1, v_2\}$ and $\overline{N}_{U_2}(v)=\{v_1\}$ for all $v \in N_U(u^*) \backslash\{v_2\}$. Then let
$$
G'=G^{\star}-\sum\limits_{v\in N_{U}(w_1)}vw_1+\sum\limits_{v\in N_{U}(w_1)}vv_1.
$$
By Lemma \ref{trans1} and $x_{v_1} \geqslant x_{w_1}$, we get $\lambda_\alpha(G')>\lambda_\alpha(G^{\star})$. Note that $G'[U\cup\{v^*\}] \cong K_b-v_1v_2$, $E_{G'}(U,W)=\{v_2w_2\}$ and $G'[W]$ is still the union of two paths at most. Hence, $G'$ is a proper subgraph of $S^{n-b}(K_b)$ and so $\lambda_\alpha(S^{n-b}(K_b))>\lambda_\alpha(G')>\lambda_\alpha(G^{\star})$, which contradicts the maximality of $\lambda_\alpha(G^{\star})$.
\end{proof}
\begin{proof}[\bf The proof of Theorem \ref{main1} {\rm (ii)}]
If $b=3$ and $\alpha\in[0,1)$, by Lemmas~\ref{lambda} and \ref{U}, we have $2=\Delta(G^{\star})\geqslant\lambda_\alpha(G^{\star})$ with equality if and only if $G^{\star}$ is a $2$-regular graph. Clearly, $C_n$ is a $2$-regular and $K_{1,3}$-minor free graph. Hence, $G^{\star}\cong C_n\cong S^{n-3}(K_3)$.

If $b\geqslant 4$ and $\alpha\in[\frac{2}{b+1},1)$, by Corollary~\ref{11} and Lemma \ref{n2v}, we have $|U_2| \geqslant|N^2(v^*)| \geqslant 2$. For $|U_2|=2$, one sees that $|N^2(v^*)|=2$. Let $U_2=\{u_1,u_2\}$. Then we have $G^{\star}[U \cup\{v^*\}] \cong K_b-u_1u_2$. By Lemma \ref{vertex}, we see that $G^{\star}[W]$ is a union of two paths at most. Suppose that $G^{\star}[W]$ consists of precisely two paths. If two pendant vertices in $W$ are connected by an edge, the resulting graph is a $K_{1, b}$-minor free graph with larger $\lambda_\alpha$-spectral radius, a contradiction. Hence, $G^{\star}[W]$ is a path with endpoints $u_1,u_2$, i.e., $G^{\star} \cong S^{n-b}(K_b)$.

For $|U_2| \geqslant 3$, let $U_{21}=\{v \in U_2|d_U(v)=b-3\}$ and $U_{22}=U_2\backslash U_{21}$. The next two claims are necessary to finish our proof.
\begin{claim}\label{u21}
If $|U_{21}|\geqslant 1$, then $G^{\star}[U_{21}]\cong K_{|U_{21}|}$. Furthermore, we have $|\overline{N}_U(v)|=|\overline{N}_{U_{22}}(v)|=1$ for all $v \in U_{21}$.
\end{claim}
\begin{proof}
If $|U_{21}|=1$, we have $G^{\star}[U_{21}]\cong K_1$. If $|U_{21}|\geqslant 2$, suppose that $G^{\star}[U_{21}]\ncong K_{|U_{21}|}$. Then there exist $v_1, v_2 \in U_{21}$ such that $v_1v_2\notin E(G^{\star})$. By the definition of $U_{21}$, one has $d_U(v_1)=d_U(v_2)=b-3$. As $|U|=b-1$, we have $N_U(v_1)\cap N_U(v_2)=U \backslash \{v_1, v_2\}$. Note that $|U_2| \geqslant 3$. Let $v_3 \in U_2 \backslash\{v_1, v_2\}$. Together with Lemma \ref{u2w}, there exists a unique vertex $w_i\in W$ such that $v_iw_i\in E(G^{\star})\ (i=1,2,3)$. Suppose there exists some $w_j\in \{w_1,w_2\}$ such that $w_j \neq w_3$. Then there is a double star with central vertices $v_j,v_3$ and $b$ leaves in $G^{\star}$, a contradiction. Thus, $w_1=w_2=w_3$.

By Lemma \ref{n2v}, $|N^2(v^*)| \geqslant 2$. Let $w_4 \in N^2(v^*) \backslash\{w_1\}$. Notice that $d_W(v)\leqslant 1$ for all $v\in U_2$ (see Lemma \ref{vertex}). Then there exists a vertex $v_4\in U_2\setminus\{v_1,v_2,v_3\}$ such that $v_4w_4\in E(G^{\star})$. Consequently, $G^{\star}$ contains a double star with central vertices $v_1,v_4$ and $b$ leaves, a contradiction. Hence, $G^{\star}[U_{21}]\cong K_{|U_{21}|}$. Then we get $\overline{N}_U(v)\subseteq U_{22}$ for all $v \in U_{21}$. Notice that $|\overline{N}_U(v)|=1$. Hence, $|\overline{N}_U(v)|=|\overline{N}_{U_{22}}(v)|=1$.
\end{proof}
\begin{claim}\label{u22-1}
If $|U_{22}| \geqslant 2$, then $u^* \in U_{21} \cup N^2(v^*)$.
\end{claim}
\begin{proof}
By Lemma \ref{ustardegree}, we have $d_{G^{\star}}(u^*)=b-1$ and there is at most one vertex in $N_{G^{\star}}(u^*)$ with a degree of at most $b-2$. Suppose that $u^* \notin U_{21} \cup N^2(v^*)$. We have {$u\in \{v^*\}\cup U_1$} or $U_{22}$ or $W\backslash N^2(v^*)$.

If $u^* \in \{v^*\}\cup U_1$, then $U_{22}\subseteq N_U(u^*)$. Note that $d_{G^{\star}}(v) \leqslant b-2$ for all $v \in U_{22}$. Then there exist at least $|U_{22}|\geqslant 2$ vertices of degree no more than $b-2$ in $N_{G^{\star}}(u^*)$, a contradiction. If $u^* \in U_{22}$, by the definition of $U_{22}$ and Lemma \ref{vertex}, we have $d_{G^{\star}}(u^*)=1+d_U(u^*)+d_W(u^*)\leqslant1+b-4+1=b-2$, a contradiction. If $u^* \in W \backslash N^2(v^*)$, by Lemma \ref{vertex}, one has $d_{G^{\star}}(u^*)=d_W(u^*)\leqslant 2$. Recall $|U_{22}| \geqslant 2$. We have $b-4\geqslant 0$ and so $d_{G^{\star}}(u^*)\leqslant b-2$, a contradiction.
\end{proof}
\begin{claim}\label{u22}
$|U_{22}| \leqslant 2$.
\end{claim}
\begin{proof}
Lemma \ref{ustardegree} states that there may be at most one vertex in $N_{G^{\star}}(u^*)$ that has a degree equal to or less than $b-2$. Suppose that $|U_{22}| \geqslant 3$. By Claim \ref{u22-1}, we have $u^* \in U_{21} \cup N^2(v^*)$.
 
If $u^* \in U_{21}$, by Claim \ref{u21}, we have $|\overline{N}_{U_{22}}(u^*)|=1$. Note that $d_{G^{\star}}(v) \leqslant b-2$ for all $v \in U_{22}$. Then there exist at least $|N_{U_{22}}(u^*)|= |U_{22}|-1\geqslant 2$ vertices in $N_{G^{\star}}(u^*)$ having degree no more than $b-2$, a contradiction.

If $u^* \in N^2(v^*)$, by $d_{G^{\star}}(u^*)=b-1$ and $d_W(u^*) \leqslant 1$ (see Lemma \ref{vertex}), we have $d_U(u^*)\geqslant b-2$. In view of Lemmas \ref{vertex} and \ref{n2v}, one sees that $d_W(v)\leqslant 1$ for all $v\in U_2$ and $|N^2(v^*)| \geqslant 2$. Therefore, $d_U(u^*)\leqslant b-2$ and so $d_U(u^*)= b-2$. By $|U|=b-1$ and $d_{G^{\star}}(v) \leqslant b-2$ for all $v \in U_{22}$, there exist at least $|U_{22}|-1\geqslant 2$ vertices in $N_{G^{\star}}(u^*)$ having degree no more than $b-2$, a contradiction.
\end{proof}

Now we come back to show Theorem \ref{main1}(ii). 

By Claim \ref{u22} and $|U_2|\geqslant 3$, we have $|U_{21}|\geqslant 1$. Furthermore, by Claim \ref{u21}, we get $G^{\star}[U_{21}]\cong K_{|U_{21}|}$ and $|\overline{N}_U(v)|=|\overline{N}_{U_{22}}(v)|=1$ for all $v \in U_{21}$. Thus, $U_{22}\neq \emptyset$ and so $1\leqslant|U_{22}| \leqslant 2$. Therefore, the following cases are being considered.

{\bf Case 1.}\ $|U_{22}| = 1$.
Let $U_{22}=\{u_0\}$ and $U_{21}=\{u_1, u_2, \ldots, u_{|U_{21}|}\}$. Clearly, $u_0u\notin E(G^{\star})$ for all $u\in U_{21}$. By Lemma \ref{u2w}, we get $d_W(u_0)=1$ and $d_W(u_i)=1\ (1\leqslant i \leqslant |U_{21}|)$. Denote by $N_W(u_i)=\{w_i\}$ where $w_i=w_j$ possibly for some $0\leqslant i\neq j \leqslant|U_{21}|$. Let $U\backslash \{u_i, u_0\}=R\ (1\leqslant i \leqslant |U_{21}|)$. Then it is easy to check that $v^* u_i \in E(G^{\star}), N_{G^{\star}}(v^*) \cap N_{G^{\star}}(u_i)=R, d_{G^{\star}}(v^*)=d_{G^{\star}}(u_i)=b-1, N_{G^{\star}}(u_i) \backslash(R \cup\{v^*\})=\{w_i\}$ and $N_{G^{\star}}(v^*) \backslash(R \cup\{u_i\})=\{u_0\}$. By Lemma \ref{trans1b} and $x_{v^*} \geqslant x_{v_i}$, we have $x_{u_0}\geqslant x_{w_i}\ (1\leqslant i \leqslant |U_{21}|)$.

Note that $G^{\star}$ is connected and $d_W(w)\leqslant 2$ for all $w\in W$ (see Lemma \ref{vertex}). Thus, every component of $G^{\star}[W]$ is a path. By Lemma \ref{n2v}, $|N^2(v^*)| \geqslant 2$. If $|N^2(v^*)|=2$, there exists a vertex $u_j \in U_{21}$ such that $N_W(u_j)=w_j\neq w_0$. Therefore, $G^{\star}$ is a subgraph of $S^{n-b}(K_b)$. If $|N^2(v^*)| \geqslant 3$, by $d_W(w)\leqslant 1$ for all $w\in N^2(v^*)$ (see Lemma \ref{vertex}), there exists a vertex $w_j\in N^2(v^*) \backslash \{w_0\}$ such that $w_j$ and $w_0$ belong to two distinct paths of $G^{\star}[W]$.
Let
$$
G'=G^{\star}-\sum_{i=1,i\neq j}^{|U_{21}|}u_iw_i+\sum_{i=1,i\neq j}^{|U_{21}|}u_iu_0.
$$
Recall $x_{u_0}\geqslant x_{w_i}$ for all $1\leqslant i \leqslant |U_{21}|$. By Lemma \ref{trans1}, $\lambda_\alpha(G')>\lambda_\alpha(G^{\star})$. Since $G'[U\cup\{v^*\}] \cong K_b-u_0u_j$, $E_{G'}(U,W)=\{u_0w_0,u_jw_j\}$ and $G'[W]$ is a union of some disjoint paths, we have $\lambda_\alpha(S^{n-b}(K_b)) \geqslant \lambda_\alpha(G')>\lambda_\alpha(G^{\star})$, a contradiction.

{\bf Case 2.}\ $|U_{22}|=2$. Let $U_{22}=\{v_1, v_2\}$. For $i=1,2$, from the definition of $U_{22}$, $d_U(v_i)\leqslant b-4$. Let $N_{U_{21}}(v_i)\supseteq\{u_i\}$. Note that $|\overline{N}_U(u_i)|=1$ and $v_1\neq v_2$. We have $u_1\neq u_2$. Consequently, $4\leqslant|U_2|\leqslant b-1$, i.e., $b\geqslant 5$. Recall $|U_{22}|=2$. By Claim \ref{u22-1}, we have $u^* \in N^2(v^*) \cup U_{21}$.

If $u^* \in U_{21}$, we may assume that $u^*=u_1$. Notice that $d_{G^{\star}}(v_2)\leqslant b-2$. Together with Lemma \ref{ustardegree} and $v_2\in N_{G^{\star}}(u_1)$, we have $d_{G^{\star}}(v_2)=b-2$ and $d_{G^{\star}}(w)=b-1$ for all $w\in N_{G^{\star}}(u_1)\setminus\{v_2\}$. By Lemma~\ref{u2w}, one has $d_W(u_1)=1$. Let $N_W(u_1)=\{w_1\}$. Then $d_{G^{\star}}(w_1)=b-1$. Notice that $w_1 \in N^2(v^*)$. Thus, $d_W(w_1) \leqslant 1$ (see Lemma~\ref{vertex}). By $|N^2(v^*)| \geqslant 2$ and $|U|=b-1$, we have $d_U(w_1)\leqslant |U|-1=b-2$. Consequently, $d_W(w_1)=b-1-d_U(w_1)\geqslant 1$. Then $d_W(w_1)=1$ and so $d_U(w_1)=b-2$. If $w_1v_1 \in E(G^{\star})$, then $G^{\star}$ contains a double star with central vertices $u_1,w_1$ and $b$ leaves, a contradiction. Hence, $w_1v_1 \notin E(G^{\star})$ and so $\overline{N}_U(w_1)=\{v_1\}$. Together with $\overline{N}_U(u_2)=\{v_2\}$. Then we see that $G^{\star}$ contains a double star with central vertices $u_2,w_1$ and $b$ leaves, a contradiction.

If $u^* \in N^2(v^*)$, by Lemma \ref{vertex}, one has $d_W(u^*) \leqslant 1$. Recall that $d(u^*)=b-1$. Then $d_U(u^*)\geqslant b-2$. By Lemma \ref{n2v}, we have $|N^2(v^*)| \geqslant 2$. Together with $d_W(v)\leqslant 1$ for all  $v\in U_2$ (see Lemma~\ref{vertex}). Then we see that $d_U(u^*)\leqslant b-2$. So, $d_U(u^*)=b-2$ and $d_W(u^*)=1$. Let $N_W(u^*)=\{w_1\}$. By Lemma~\ref{ustardegree} and $b \geqslant 5$, we have $d_{G^{\star}}(w_1) \geqslant b-2 \geqslant 3$. Suppose that $w_1 \notin N^2(v^*)$. By Lemma~\ref{vertex}, one has $d_{G^{\star}}(w_1)=d_W(w_1) \leqslant 2$, a contradiction. Therefore, $w_1 \in N^2(v^*)$. By Lemma~\ref{vertex} and $d_U(u^*)=b-2$, we have $d_W(w_1) \leqslant 1$ and $d_U(w_1)=1$. So, $d_{G^{\star}}(w_1)\leqslant 2$, a contradiction.
\end{proof}
\section{\normalsize  The proof of Theorem~\ref{main2}.}
In this section, we give the proof of Theorem~\ref{main2}, which characterizes the graph $G^*$ having the maximum $A_\alpha$-spectral radius among all $n$-vertex $K_{a,b}$-minor free graphs with $\alpha\in[0,1)$, $a\geqslant 2$ and sufficiently large $n$. The following lemma is derived from the arguments of \cite[Theorem~1.2]{Zhang} and \cite[Theorem~1.2]{Tait}.
\begin{lem}[\cite{Zhang, Tait}]\label{S}
$G^*$ has a clique dominating set $S^*$ of cardinality $a-1$.
\end{lem}
First, we give an important property of $G^*-S^*$.
\begin{lem}\label{emax}
Assume that $\mathcal{F}$ is made up of a few components in $G^*-S^*$ with $|V(\mathcal{F})|<c_1$, where $c_1$ is a constant. If $\mathcal{F'}$  is a graph with $(a,b)$-property satisfying $V(\mathcal{F'})=V(\mathcal{F}),$  then $e(\mathcal{F'})\leqslant e(\mathcal{F})$.
\end{lem}
\begin{proof}
Let $G'=G^*-E(\mathcal{F})+E(\mathcal{F'})$. Notice that $X_M=\max_{v\in V(\mathcal{F})}x_v,\ X_m=\min_{v\in V(\mathcal{F})}x_v$. Then 
\begin{align*}
\lambda_\alpha(G')-\lambda_\alpha(G^*)
& \geqslant {\bf x}^T\left(A_\alpha(G')-A_\alpha(G^*)\right){\bf x}\\
&=\sum_{u v \in E(\mathcal{F'})}(\alpha(x_u^2+x_v^2)+2(1-\alpha)x_ux_v)
-\sum_{u v \in E(\mathcal{F})}(\alpha(x_u^2+x_v^2)+2(1-\alpha)x_ux_v)\\
& \geqslant 2 e(\mathcal{F'})X_m^2-2 e(\mathcal{F})X_M^2.
\end{align*}
As $G^*$ is $K_{a, b}$-minor free, together with Lemma~\ref{ab}, we see that $\mathcal{F}$ has the $(a,b)$-property. Thus, $\mathcal{F}$ is $K_{1,b}$-minor free and so $\Delta(\mathcal{F})<b$. By Lemma~\ref{S}, we have $|S^*|=a-1$ is a constant, one sees that $G^*$ and $\mathcal{F}$ satisfy the conditions of Lemma \ref{Mm}. Suppose that $e(\mathcal{F'})>e(\mathcal{F})$. By Lemma \ref{Mm}, $\lambda_\alpha(G')-\lambda_\alpha(G^*)\geqslant 2 e(\mathcal{F'})X_m^2-2 e(\mathcal{F})X_M^2>0.$ As $\mathcal{F'}$ has the $(a,b)$-property, together with~Lemma~\ref{ab}, one sees that $G'$ is also $K_{a,b}$-minor free. Then we get a contradiction to the maximality of $\lambda_\alpha(G^*)$. So, $e(\mathcal{F'})\leqslant e(\mathcal{F})$.
\end{proof}
Let $\mathscr{F}_{i},\ \mathscr{F}_{\neq i}$ and $\mathscr{F}_{>i}$ represent the sets of components in $G^*-S^*$  with orders equal to $i$, not equal to $i$, and larger than $i$, respectively. Next, we give two lemmas to determine the number of some special family of components in $G^*-S^*$.
\begin{lem}\label{b23}
If $b\leqslant 3$, then $|\mathscr{F}_{\neq b}|\leqslant 1$. In particular, if $b=2$, we have $\mathscr{F}_{\neq2}=\mathscr{F}_1$.
\end{lem}
\begin{proof}
Let $F$ be one component in $G^*-S^*$. By Lemma \ref{ab}, we know that $F$ has the $(a,b)$-property. Thus, $F$ is $K_{1,b}$-minor free and $\Delta(F)\leqslant b-1$. If $b=2$, we get $\Delta(F)\leqslant 1$. Then $F\cong K_2$ or $K_1$. So, $|\mathscr{F}_{>2}|=0$ and $\mathscr{F}_{\neq2}=\mathscr{F}_1$. If $b=3$, let $|V(F)|=h$. By $\Delta(F)\leqslant 2$, we have $F\cong C_h$ or $P_h$. Notice that $G^*=K_{a-1}\vee F$ is $K_{a,b}$-minor free. Then $F\cong C_3$ or $P_h$.

Suppose that there exist two components $F_1, F_2 \in \mathscr{F}_{\neq b}$. We choose $u\in V(F_1)$ and $v\in V(F_2)$ such that $d_{F_1}(u)=d_{F_2}(v)=1$. Then $F_1\cup F_2+uv$ also has the $(a,b)$-property. By Lemma \ref{ab}, $G^*+uv$ is also $K_{a,b}$-minor free. Together with Lemma \ref{sub}, we have $\lambda_\alpha(G^*+uv)>\lambda_\alpha(G^*)$, which contradicts the maximality of $\lambda_\alpha(G^*)$. Hence, $|\mathscr{F}_{\neq b}|\leqslant 1$.
\end{proof}
\begin{lem}\label{b4>b}
If $b\geqslant 4$, then we have
\begin{wst}
\item[{\rm (i)}]$\mathscr{F}_{>b+3}=\emptyset$;
\item[{\rm (ii)}]$|\mathscr{F}_{<b}\cup\mathscr{F}_{b+2}|\leqslant1$;
\item[{\rm (iii)}]$|\mathscr{F}_b|=O\left(\frac{n}{b}\right)$.
\end{wst}
\end{lem}
\begin{proof}
{\rm (i)}\ Suppose that there exists a component $F\in \mathscr{F}_{>b+3}$ with order $pb+q\ (0\leqslant q\leqslant b-1)$. As $F$ is $K_{1,b}$-minor free, together with Lemma~\ref{edge1b}, we have $e(F)\leqslant {b\choose 2}+(pb+q)-b$. Let $F'=pK_b\cup K_q$ such that $V(F')=V(F)$. Clearly, $F'$ has the $(a,b)$-property and $e(F')=p{b\choose 2}+{q\choose 2}$. Note that $pb+q> b+3$ with $b\geqslant 4$. Then we have
\begin{align*}
e(F')-e(F)\geqslant&p{b\choose 2}+{q\choose 2}-{b\choose 2}-(pb+q)+b\\
=&\frac{(p-1)b(b-3)+q(q-3)}{2}\\
\geqslant& \frac{4(p-1)+q(q-3)}{2}>0.
\end{align*}
So, $e(F')>e(F)$.

If $p\leqslant 3$, then $|V(F)|=|V(F')|<4b$, where $4b$ is a constant. By Lemma \ref{emax}, $e(F')\leqslant e(F)$, a contradiction. If $p\geqslant 4$, by $b\geqslant 4$, we have
$$
\frac{11}{12}e(F')-e(F)\geqslant \frac{11}{12}p{b\choose 2}-{b\choose 2}-pb
=\frac{b((11p-12)b-35p+12)}{24}\geqslant \frac{9b(p-4)}{24}\geqslant 0.
$$
Hence, $\frac{11}{12}e(F')\geqslant e(F)$. Let $G'=G^*-E(F)+E(F')$. By Lemma \ref{Mm}, we have
\begin{align*}
\lambda_\alpha(G')-\lambda_\alpha(G^*)
\geqslant& {\bf x}^T(A_\alpha(G')-A_\alpha(G^*)){\bf x}\\
=&\sum_{uv\in E(F')}\left(\alpha(x_u^2+x_v^2)+2(1-\alpha)x_ux_v\right)
-\sum_{uv\in E(F)}\left(\alpha(x_u^2+x_v^2)+2(1-\alpha)x_ux_v\right)\\
\geqslant& 2e(F')X_m^2-2e(F)X_M^2\\
\geqslant& 2e(F')\left(X_m^2-\frac{11}{12}X_M^2\right)>0.
\end{align*}
By Lemma \ref{ab} and $F'$ with the $(a,b)$-property, $G'$ is $K_{a,b}$-minor free.
Then we get a contradiction. Thus, $\mathscr{F}_{>t+3}=\emptyset$.

{\rm (ii)}\ By Lemmas \ref{edge1b} and \ref{emax}, for all components $F\in \mathscr{F}_{<b}\cup\mathscr{F}_{b+2}$, we have $e(F)\leqslant {|V(F)|\choose 2}$ if $F\in \mathscr{F}_{<b}$ and $e(F)={b\choose 2}+2$ otherwise. Suppose that $|\mathscr{F}_{<b}\cup\mathscr{F}_{b+2}|\geqslant 2$. We choose $F_1, F_2 \in \mathscr{F}_{<b}\cup\mathscr{F}_{b+2}$ such that $|V(F_1)|+|V(F_2)|=pb+q\ (0 \leqslant p \leqslant 2,\ 0 \leqslant q\leqslant b-1)$. It is easy to get that $e(F_1 \cup F_2)<p{b\choose 2}+{q\choose 2}=e(pK_b \cup K_q)$. Note that $pK_b \cup K_q$ has the $(a,b)$-property and $|V(F_1 \cup F_2)|=|V(pK_b \cup K_q)|$. We get a contradiction to Lemma \ref{emax}.

{\rm (iii)}\
By {\rm (i)} and {\rm (ii)}, we just need to show that $\left|\mathscr{F}_{b+1}\right|<b$ and $\left|\mathscr{F}_{b+3}\right|<b$. Suppose that there exists $\left|\mathscr{F}_j\right| \geqslant b$ for some $j\in \{b+1,\,b+3\}$. Let $\mathcal{F}$ represent a set that consists of disjoint union of $b$ components in $\mathscr{F}_j$. As $jK_b$ has the $(a,b)$-property, together with Lemma \ref{emax}, we have $e(jK_b)\leqslant e(\mathcal{F})$. To complete our proof, it suffices to show that $e(\mathcal{F})<e(jK_b)$.

If $j=b+1$, by Lemma \ref{eb+1} and $\tau=\left\lfloor\frac{b+1}{a+1}\right\rfloor<\frac{b+1}{2}$, we have
$$
e(\mathcal{F})\leqslant b\left({b\choose 2}+\tau-1\right)<(b+1){b\choose 2}=e(jK_b).
$$

If $j =b+3$, by Lemma \ref{edge1b}, one has
$$
e(\mathcal{F}) \leqslant b\left({b\choose 2}+3\right)<(b+3){b\choose 2} =e(jK_b).
$$

This completes the proof.
\end{proof}

\begin{lem}\label{edgeb23}
If $b\geqslant 4$, we have $e(F)>{b\choose 2}+3$ for any component $F\in \mathscr{F}_{b+3}$.
\end{lem}
\begin{proof}
As $F$ has the $(a,b)$-property, together with Lemma \ref{emax},  we have $e(F)\geqslant e(K_b\cup K_3)={b\choose 2}+3$. Suppose that $e(F)\leqslant {b\choose 2}+3$. Then $e(F)={b\choose 2}+3$.

By Lemma \ref{b4>b}, we have $F\cup 2K_b\subseteq G^*-S^*$. Suppose that $\tau\geqslant 3$. Since $e(\overline{F_{a,b}})={b+1\choose 2}-e(F_{a,b})={b\choose 2}+\tau-1$, one has
$$
e(F\cup 2K_b)=3{b\choose 2}+3<3{b\choose 2}+3(\tau-1)=e(3\overline{F_{a,b}}).
$$
Note that $|V(\overline{F_{a,b}})|=b+1$ and $\overline{F_{a,b}}$ has the $(a,b)$-property (see Lemma \ref{abgraph}). By Lemma \ref{emax}, $e(3\overline{F_{a,b}})\leqslant e(F\cup 2K_b)$,
a contradiction. Thus, $\tau\leqslant 2$.

Suppose there exists a vertex $v\in V(F)$ such that $d_{F}(v)=1$. Clearly, $F-v$ also has the $(a,b)$-property and $e(F-v)={b\choose 2}+2$. By Lemma \ref{b+21} and $\tau\leqslant 2$, we have $F-v$ is isomorphic to either some $F(a_1,a_2,a_3)$ or $\overline{P^{\star}}$. If $F-v$ is isomorphic to some $F(a_1,a_2,a_3)$, by $\Delta(F)<b$, then $v$ is adjacent to some $v_i\ (i=1,2,3)$. Contradicting the edges $v_1v_2$ and $v_2v_3$ gives us a graph, which contains a copy of $K_{1,b}$, a contradiction. If $F-v$ is isomorphic to $\overline{P^{\star}}$, assume that $F$ contains a path $vuw$. Obviously, $u$ and $w$ are non-adjacent vertices in $P^{\star}$ and $|N_{P^{\star}}(u)\cap N_{P^{\star}}(w)|=1$. Recall that $|V(F)|=b+3$. We have $|N_{F}(u)\cup N_{F}(w)|=b+2$. Consequently, $F$ contains a double star with central vertices $u,w$ and $b$ leaves, a contradiction. Therefore, $d_{F}(v)\geqslant 2$ for all $v\in V(F)$. For convenience, let $F'=K_b\cup K_3$.

Denote by $d_1, d_2, \ldots, d_{b+3}$ and $d'_1, d'_2, \ldots, d'_{b+3}$ the degree sequences of $F$ and $F'$, respectively, in decreasing order. Let $X^T=\{d_1, d_2, \ldots, d_{b+3})$ and $Y^T=(d'_1, d'_2, \ldots, d'_{b+3})$. For the sake of simplicity and without making any assumptions, let $V(F)=V(F')$. Then we have
\begin{equation}\label{1}
\begin{array}{ll}
\sum\limits_{i=1}^{b+3}d_id_i
&=\sum\limits_{v\in V(F)}d_F(v)d_{F}(v)
=2\sum\limits_{uv\in E(F)}d_F(v)
=\sum\limits_{uv\in E(F)}(d_F(u)+d_F(v)),\\
\sum\limits_{i=1}^{b+3}d_id'_i
&=\sum\limits_{v\in V(F')}d_F(v)d_{F'}(v)
=2\sum\limits_{uv\in E(F')}d_F(v)
=\sum\limits_{uv\in E(F')}(d_F(u)+d_F(v))\\
&=\sum\limits_{v\in V(F)}d_F(v)d_{F'}(v)
=2\sum\limits_{uv\in E(F)}d_{F'}(v)
=\sum\limits_{uv\in E(F)}(d_{F'}(u)+d_{F'}(v)),\\
\sum\limits_{i=1}^{b+3}d'_id'_i
&=\sum\limits_{v\in V(F')}d_{F'}(v)d_{F'}(v)
=2\sum\limits_{uv\in E(F')}d_{F'}(v)
=\sum\limits_{uv\in E(F')}(d_{F'}(u)+d_{F'}(v)).
\end{array}
\end{equation}
By $e(F)=e(F')$, we have $\sum_{i=1}^{b+3}d_i=\sum_{i=1}^{b+3}d'_i$. Together with $2\leqslant d_i\leqslant b-1$ for $i=1,2,\ldots,b+3$, it is easy to get $X \prec Y$. By Lemma \ref{major2}, we have
$$
\sum\limits_{i=1}^{b+3}d_id_i=X^T \cdot X \leqslant Y^T \cdot X=\sum\limits_{i=1}^{b+3}d'_id_i.
$$
Let $G'=G^*-E(F)+E(F')$. Clearly, $G'$ is $K_{a,b}$-minor free. We proceed by considering the following cases.

{\bf Case 1.}\ $X\neq Y$. Let $d_B=\frac{1}{b} \sum\limits_{i=1}^{b}d_i$ and $d_A=\frac{1}{3} \sum\limits_{i=b+1}^{b+3}d_i$. Combing with $2\leqslant d_i\leqslant b-1$ and $X\neq Y$ gives us $d_A>2$. Furthermore, we have
$$
\left\{
\begin{array}{ll}
2e(F)=bd_B+3d_A=b(b-1)+6,\\
\sum\limits_{i=1}^{b+3}d_i'd_i'=3\times 2^2+b(b-1)^2,\\
\sum\limits_{i=1}^{b+3}d_id_i'=2\times 3d_A+b(b-1)d_B.
\end{array}
\right.
$$
Together with $b\geqslant 4$, we have $\sum\limits_{i=1}^{b+3}d_i'd_i'-\sum\limits_{i=1}^{b+3}d_id_i'=3(b-3)(d_A-2)>0$. Consequently,
\begin{equation}\label{2}
\sum\limits_{i=1}^{b+3}d_id_i\leqslant \sum\limits_{i=1}^{b+3}d_id'_i\leqslant\sum\limits_{i=1}^{b+3}d'_id'_i-1.
\end{equation}

Let $\lambda_\alpha(G'){\bf y}=A_\alpha(G'){\bf y}$ and denote by $Y_s=\sum\limits_{v\in S^*}y_s$, $Y_M=\max_{v\in V(F')}y_v$ and $Y_m=\min_{v\in V(F')}y_v.$ By \eqref{F} and \eqref{X}, for all $v\in  V(F)=V(F')$, one has
{\footnotesize
$$\frac{(1-\alpha)(X_s+d_F(v)X_m)}{\lambda_\alpha(G^*)-\alpha (a+1)}
\leqslant\frac{(1-\alpha)(X_s+d_F(v)X_m)}{\lambda_\alpha(G^*)-\alpha (a-1+d_F(v))}
\leqslant x_v
\leqslant \frac{(1-\alpha)(X_s+d_F(v)X_M)}{\lambda_\alpha(G^*)-\alpha (a-1+d_F(v))}
\leqslant \frac{(1-\alpha)(X_s+d_F(v)X_M)}{\lambda_\alpha(G^*)-\alpha (a+b-2)}. 
$$
$$\frac{(1-\alpha)(Y_S+d_{F'}(v)Y_m)}{\lambda_\alpha(G')-\alpha (a+1)}
\leqslant\frac{(1-\alpha)(Y_S+d_{F'}(v)Y_m)}{\lambda_\alpha(G')-\alpha (a-1+d_{F'}(v))}
\leqslant y_v
\leqslant \frac{(1-\alpha)(Y_S+d_{F'}(v)Y_M)}{\lambda_\alpha(G')-\alpha (a-1+d_{F'}(v))}
\leqslant \frac{(1-\alpha)(Y_S+d_{F'}(v)Y_M)}{\lambda_\alpha(G')-\alpha (a+b-2)},
$$
}
and
{\footnotesize
\begin{align}\label{xyMm}
&X_m\geqslant\frac{(1-\alpha)X_s}{{ \lambda_\alpha(G^*)}-\alpha(a-1)};\ X_M<\frac{(1-\alpha)X_s}{{ \lambda_\alpha(G^*)}-\alpha(a-1)-b};\ Y_m\geqslant\frac{(1-\alpha)Y_s}{{ \lambda_\alpha(G')}-\alpha(a-1)};\ Y_M<\frac{(1-\alpha)Y_s}{{\lambda_\alpha(G')}-\alpha(a-1)-b}.
\end{align}
}
By some calculations, we have
\begin{align*}
{\bf x}^{T} {\bf y}(\lambda_\alpha(G')-\lambda_\alpha(G^*))=&{\bf x}^{T}(A_{\alpha}(G')-A_{\alpha}(G^*)) {\bf y}\\
= &\sum\limits_{uv \in E(F')}\big(\alpha(x_uy_u+x_vy_v)+(1-\alpha)(x_uy_v+x_vy_u)\big)\\
&-\sum\limits_{uv \in E(F)}\big(\alpha(x_uy_u+x_vy_v)+(1-\alpha)(x_uy_v+x_vy_u)\big) \\
= &\alpha \big(\sum\limits_{uv \in E(F')}(x_uy_u+x_vy_v)-\sum\limits_{uv \in E(F)}(x_uy_u+x_vy_v)\big)\\
&+(1-\alpha)\big(\sum\limits_{uv \in E(F')}(x_uy_v+x_vy_u)-\sum\limits_{uv \in E(F)}(x_uy_v+x_vy_u)\big) \\
\geqslant&\left(\frac{\alpha { (1-\alpha)^2}B_1}
{(\lambda_\alpha(G^*)-\alpha (a+1))(\lambda_\alpha(G')-\alpha (a+1))}\right.\\
&\left.-\frac{\alpha{ (1-\alpha)^2} B_2}{(\lambda_\alpha(G^*)-\alpha (a+b-2))(\lambda_\alpha(G')-\alpha (a+b-2){\b)}}\right)\\
&+\left(\frac{{ (1-\alpha)^3}B_3}{(\lambda_\alpha(G^*)-\alpha (a+1))(\lambda_\alpha(G')-\alpha (a+1)})\right.\\
&-\left.\frac{{ (1-\alpha)^3}B_4}{(\lambda_\alpha(G^*)-\alpha (a+b-2))(\lambda_\alpha(G')-\alpha (a+b-2){)}}\right),
\end{align*}
where
\begin{align*}
B_1=&{ \sum\limits_{uv\in E(F')}((X_s+d_{F}(u)X_m)(Y_s+d_{F'}(u)Y_m)+(X_s+d_{F}(v)X_m)(Y_s+d_{F'}(v)Y_m))},\\
B_2=&{ \sum\limits_{uv\in E(F)}((X_s+d_{F}(u)X_M)(Y_s+d_{F'}(u)Y_M)+(X_s+d_{F}(v)X_M)(Y_s+d_{F'}(v)Y_M)),}\\
B_3=&\sum\limits_{uv \in E(F')}\left((X_s+d_{F}(u)X_m)(Y_s+d_{F'}(v)Y_m)+(X_s+d_{F}(v)X_m)(Y_s+d_{F'}(u)Y_m)\right),\\
B_4=&\sum\limits_{uv \in E(F)}\left((X_s+d_{F}(u)X_M)(Y_s+d_{F'}(v)Y_M)+(X_s+d_{F}(v)X_M)(Y_s+d_{F'}(u)Y_M)\right).
\end{align*}
By some calculations, we see
\begin{eqnarray*}
  B_1 &\geqslant&  2e(F')X_sY_s+\sum\limits_{i=1}^{b+3}d_id'_i X_mY_s+\sum\limits_{i=1}^{b+3}d'_id'_iX_sY_m \\
  B_2 &\leqslant& 2e(F)X_sY_s+ \sum\limits_{i=1}^{b+3}d_id_iX_MY_s +\sum\limits_{i=1}^{b+3}d_id'_iX_sY_M+c_2X_MY_M\\
      &\leqslant& 2e(F)X_sY_s+ \sum\limits_{i=1}^{b+3}d_id_i'X_MY_s +\sum\limits_{i=1}^{b+3}d_id'_iX_sY_M+c_2X_MY_M,\\
  B_3 &\geqslant& 2e(F')X_sY_s+ \sum\limits_{i=1}^{b+3}d_id'_i X_mY_s + \sum\limits_{i=1}^{b+3}d'_id'_i X_sY_m,  \\
  B_4 &\leqslant& 2e(F)X_sY_s+\sum\limits_{i=1}^{b+3}d_id'_i X_MY_s + \sum\limits_{i=1}^{b+3}d_id'_i X_sY_M+c_3X_MY_M,
\end{eqnarray*}
where $c_2=\sum_{uv\in E(F)}(d_{F}(u)d_{F'}(u)+d_{F}(v)d_{F'}(v))$ and $c_3=\sum_{uv\in E(F)}d_F(u)d_{F'}(v)+d_F(v)d_{F'}(u)$ are constants.
Together with Lemma \ref{Mm}, $e(F)=e(F')$, \eqref{1}, \eqref{2} and \eqref{xyMm}, we have
\begin{align*}
B_1-B_2
\geqslant& \left(\sum\limits_{i=1}^{b+3}d'_id'_i-\frac{1}{2}\right)X_sY_m-\sum\limits_{i=1}^{b+3}d_id'_iX_sY_M
+\frac{1}{2}X_sY_m+\sum\limits_{i=1}^{b+3}d_id'_iY_s(X_m-X_M)-c_2X_MY_M\\
>&\frac{(1-\alpha)X_sY_s}{2}\left(\frac{1}{\lambda_\alpha(G')-\alpha(a-1)}
+\frac{2\sum\limits_{i=1}^{b+3}d_id'_i}{\lambda_\alpha(G^*)-\alpha(a-1)}
-\frac{2\sum\limits_{i=1}^{b+3}d_id'_i}{\lambda_\alpha(G^*)-\alpha(a-1)-b}\right.\\
&\left.-\frac{c_2}{(\lambda_\alpha(G^*)-\alpha(a-1)-b)(\lambda_\alpha(G')-\alpha(a-1)-b)} \right)
>0,
\end{align*}
i.e., $B_1>B_2$. Similarly, we get $B_3-B_4>0$. Note that $\lambda_\alpha(G)\geqslant \sqrt{n-1}$, $\alpha, a, b$ are all constants and $n$ is large enough. Then we have
\begin{align*}
\frac{B_1}{(\lambda_\alpha(G^*)-\alpha (a+1))(\lambda_\alpha(G')-\alpha (a+1))}
-\frac{B_2}{(\lambda_\alpha(G^*)-\alpha (a+b-2))(\lambda_\alpha(G')-\alpha (a+b-2))}>0
\end{align*}
and
\begin{align*}
\frac{B_3}{(\lambda_\alpha(G^*)-\alpha (a+1))(\lambda_\alpha(G')-\alpha (a+1))}
-\frac{B_4}{(\lambda_\alpha(G^*)-\alpha (a+b-2))(\lambda_\alpha(G')-\alpha (a+b-2))}>0.
\end{align*}
Consequently, ${\bf x}^{T} {\bf y}(\lambda_\alpha(G')-\lambda_\alpha(G^*))>0.$ Recall $G'$ is $K_{a,b}$-minor free. Then we get a contradiction to the maximality of $\lambda_\alpha(G^*)$.

{\bf Case 2.}\ $X=Y$. Let $A=\{v\in V(F)|d_F(v)=2\}$, $B_1=N_F(A)$ and $B_2=V(F)\setminus (A\cup B_1)$. Clearly, $|B_1\cup B_2|=b$ and $d_F(u)=b-1$ for all $u\in B_1\cup B_2$. Moreover, we have $d_A(u)=|\overline{N}_{B_1\cup B_2}(u)|$ if $u\in B_1$ and $N_F(u)\cap \{u\}=B_1\cup B_2$ if $u\in B_2$.
Since $F$ is connected, we have $2\leqslant e(A,B_1)\leqslant 6$, i.e., there exist at most three non-adjacent vertex-pair in $B_1\cup B_2$.

Suppose that there exists a vertex $u\in B_1$ such that $d_A(u)\geqslant 2$. If $B_2\neq \emptyset$, let $w\in B_2$. Since $N_F(w)\cap \{w\}=B_1\cup B_2$, we have $wu\in E(F)$. Consequently, $G^*$ contains a double star with central vertices $w,u$ and $b$ leaves, a contradiction. If $B_2=\emptyset$ and $d_A(u)=3$, let $\overline{N}_{B_1}(u)=\{u_1,u_2,u_3\}$. Recall that there are at most three non-adjacent vertex-pair in $B_1$. We have $\overline{N}_{B_1}(u_i)=\{u\}$ and so $d_A(u_i)=1\ (i=1,2,3)$. Assume that $N_A(u_1)=\{v_1\}$. Contracting the edge $v_1u_1$ to a vertex $v$ gives us the graph containing a double star with central vertices $u,v$ and $b$ leaves, a contradiction. If $B_2=\emptyset$ and $d_A(u)=2$, by $b\geqslant 3$, let $\overline{N}_{B_1}(u)=\{u_1,u_2\}$ and $u_3\in N_{B_1}(u)$. Suppose that $u_1,u_2\in N_F(u_3)$. We obtain that $F$ contains a double star with central vertices $u,u_3$ and $b$ leaves, a contradiction. Recall that there are at most three non-adjacent vertex pairs in $B_1$. Thus, we see that either $u_1u_3\notin E(F)$ or $u_2u_3\notin E(F)$. We may consider $u_1u_3\notin E(F)$. Then we have $\overline{N}_{B_1}(u_1)=\{u, u_3\}$ and so $d_A(u_1)=2$. Note that $|A|=3$. Thus, $N_A(u)\cap N_A(u_1)\neq \emptyset$. If $N_A(u)=N_A(u_1)$, then contract the edge $u_2u_3$. Consequently, $K_{2,b-1}$, with one partite set $\{u,v\}$, is a subgraph of the resultant graph, a contradiction. If $N_A(u)\neq N_A(u_1)$, let $v_1\in N_A(u)\cap N_A(u_1)$. Contracting the edge $v_1u_1$ to a vertex, say $v,$ gives us a graph containing a double star with central vertices $u,v$ and $b$ leaves, a contradiction. Therefore, we have $d_A(u)=1$ for all $u\in B_1$. So, for each $u_i \in B_1$, we have $|\overline{N}_{B_1}(u_i)|=1$ and assume that $N_{A}(u_i)=\{v_i\}$.

Assume that $u_i,u_j\in B_1$ and $u_iu_j\notin E(F)$. If $v_i=v_j$, then $F$ contains a copy of $K_{2,b-1}$, a contradiction. If $v_i \neq v_j$ and $v_i v_j \in E(H)$, by contracting the edge $v_i v_j$, $K_{2,b-1}$ is a subgraph of the resulting graph, a contradiction. Then for each $u_i,u_j\in B_1$ and $u_iu_j\notin E(F)$, $v_i \neq v_j$ and $v_i v_j \notin E(H)$. Moreover
\begin{align*}
\lambda_\alpha(G')-\lambda_\alpha(G^*)\geqslant {\bf x}^T\left(A_\alpha(G')-A_\alpha(G^*)\right){\bf x}
=2(1-\alpha)\sum\limits_{u_iu_j\notin E(H^{\star})}(x_{u_i}-x_{v_j})(x_{u_j}-x_{v_i}).
\end{align*}
By Lemma \ref{Mm} and $d_{F}(u_i)=d_{F}(u_j)=b-1>2=d_{F}(v_j)=d_{F}(v_i)$, one has $x_{u_i}>x_{v_j}$ and $x_{u_j}>x_{v_i}$. From Lemma \ref{trans1}, we obtain $\lambda_\alpha( G')>\lambda_\alpha(G^*)$, a contradiction.
\end{proof}
\begin{lem}\label{>b+2}
If $b\geqslant 4$, then $\mathscr{F}_{>b+2}=\emptyset$.
\end{lem}
\begin{proof}
Suppose there exists a component $F\in \mathscr{F}_{b+3}$. By $F$ is $K_{1,b}$-minor free and Lemma \ref{edge1b}, $e(F)\leqslant {b\choose 2}+3$. On the other hand, by Lemma \ref{edgeb23}, we get $e(F)>{b\choose 2}+3$, a contradiction. Thus, $\mathscr{F}_{b+3}=\emptyset$. Together with $\mathscr{F}_{>b+3}\neq\emptyset$ (see Lemma \ref{b4>b}), we have $\mathscr{F}_{>b+2}=\emptyset$.
\end{proof}
Next, we consider the components in $\mathscr{F}_{b+2}$ and $\mathscr{F}_{b+1}$, respectively.
\begin{lem}\label{b+2}
If $b\geqslant 4$ and $\mathscr{F}_{b+2}\neq \emptyset$, we have $\tau=\left\lfloor\frac{b+1}{a+1}\right\rfloor\leqslant 2$. Furthermore, for any component $F\in \mathscr{F}_{b+2}$, we have $F$ is isomorphic to $F(a,0,b-1-a)$ if $\tau =2$ and $\overline{P^{\star}}$ otherwise.
\end{lem}
\begin{proof}
Note that $F\in \mathscr{F}_{b+2}$ is $K_{1,b}$-minor free. In view of Lemma \ref{edge1b}, $e(F)\leqslant {b\choose 2}+2$. Note that $|V(F(\omega, 0,b-1-\omega))|=b+2$ and $F(\omega, 0,b-1-\omega)$ has the $(a,b)$-property (see Lemma \ref{abgraph}). Together with Lemma \ref{emax}, $e(F)\geqslant e(F(\omega, 0,b-1-\omega))={b\choose 2}+2$. Thus, $e(F)={b\choose 2}+2$.

By Lemma~\ref{b4>b} (iii), $K_b\cup F\subseteq G^*-S^*$. Let $F'=F(a,b)\cup F(a,b)$. Since $F(a,b)$ has the $(a,b)$-property (see Lemma \ref{abgraph}), one sees that $F'$ also has the $(a,b)$-property. Together with Lemma \ref{emax} and $|V(F')|=|V(K_b\cup F)|$, we have $e(F')\leqslant e(K_b\cup F)$. Suppose that $\tau \geqslant 3$. By Lemma \ref{eb+1}, $e(F(a,b))={b\choose 2}+\tau-1 \geqslant {b\choose 2}+2$. Thus, $e(K_b\cup F)=2{b\choose 2}+2<e(F'),$ a contradiction. Therefore, $\tau \leqslant 2$.
\unitlength 1.7mm
\begin{center}
\begin{picture}(90,50)\linethickness{0.8pt}
\put(60,12){$F(a_1,a_2,a_3)$}

\put(42,42){$v_1$}
\put(42.5,41){\circle*{1}}
\Line(42.5,41)(44.5,31.1)
\Line(42.5,41)(43.2,30.9)

\put(63.5,42){$v_2$}
\put(64,41){\circle*{1}}
\Line(64,41)(64.4,30.4)
\Line(64,41)(62.9,30.4)

\put(84.5,42){$v_3$}
\put(85,41){\circle*{1}}
\Line(85,41)(82.7,30.7)
\Line(85,41)(84,30.7)
\Line(42.5,41)(64,41)
\Line(64,41)(85,41)

\put(40,20.5){$K_{b-1}$}
\polygon(40,19.7)(89,19.7)(89,33.4)(40,33.4)

\put(43,26.6){$V_1$}
\cbezier(40.5,27.3)(40.8,20)(49.8,20)(50.1,27.3)\cbezier(40.5,27.3)(40.8,34.5)(49.8,34.5)(50.1,27.3)
\put(44.5,31.1){\circle*{.8}}
\put(43.2,30.9){\circle*{.8}}

\put(62,26.1){$V_2$}
\cbezier(59,27.3)(59.3,20)(68.2,20)(68.5,27.3)\cbezier(59,27.3)(59.3,34.5)(68.2,34.5)(68.5,27.3)
\put(64.4,30.4){\circle*{.8}}
\put(62.9,30.4){\circle*{.8}}

\put(83,25.6){$V_3$}
\cbezier(78.2,27)(78.5,19.7)(87.4,19.7)(87.7,27)\cbezier(78.2,27)(78.5,34.2)(87.4,34.2)(87.7,27)
\put(82.7,30.7){\circle*{.8}}
\put(84,30.7){\circle*{.8}}

\put(46.7,23.7){\circle*{.8}}
\put(47.5,24.5){\circle*{.8}}
\cbezier(46.7,23.7)(57.9,19.8)(70.3,19.2)(81.3,24.1)
\cbezier(47.5,24.5)(58.3,20)(70.6,20.6)(80.8,25.3)
\put(80.8,25.3){\circle*{.8}}
\put(81.3,24.1){\circle*{.8}}
\put(48.2,27.3){\circle*{.8}}
\put(48.2,26.4){\circle*{.8}}
\Line(48.2,27.3)(59.6,27.3)
\Line(48.2,26.4)(59.8,26.3)
\put(59.6,27.3){\circle*{.8}}
\put(59.8,26.3){\circle*{.8}}
\put(67.8,26.8){\circle*{.8}}
\put(67.7,27.7){\circle*{.8}}
\Line(67.8,26.8)(79.4,26.7)
\Line(67.7,27.7)(79.2,27.7)
\put(79.4,26.7){\circle*{.8}}
\put(79.2,27.7){\circle*{.8}}

\put(15,12){$P^{\star}$}

\put(-1,34.3){$v_5$}
\put(1.7,35.2){\circle*{1}}
\Line(1.7,35.2)(7.1,21.2)
\Line(1.7,35.2)(6.7,35.3)

\put(6,33){$u_5$}
\put(6.7,35.3){\circle*{1}}
\Line(6.7,35.3)(24,35.4)

\put(7,19.1){$v_4$}
\put(7.1,21.2){\circle*{1}}
\Line(7.1,21.2)(23.9,21)

\put(11,23.8){$u_4$}
\put(9.8,24.9){\circle*{1}}
\Line(9.8,24.9)(15.6,40.2)
\Line(9.8,24.9)(7.1,21.2)

\put(16.5,39.5){$u_1$}
\put(15.6,40.2){\circle*{1}}
\Line(15.6,40.2)(21.4,24.8)

\put(16.5,45.7){$v_1$}
\put(15.6,46){\circle*{1}}
\Line(15.6,46)(15.6,40.2)
\Line(15.6,46)(29.1,35.4)
\Line(15.6,46)(1.7,35.2)

\put(22.5,24.5){$u_3$}
\put(21.4,24.8){\circle*{1}}
\Line(21.4,24.8)(23.9,21)
\Line(21.4,24.8)(6.7,35.3)

\put(24,19){$v_3$}
\put(23.9,21){\circle*{1}}

\put(24,33.4){$u_2$}
\put(24,35.4){\circle*{1}}
\Line(24,35.4)(29.1,35.4)
\Line(24,35.4)(9.8,24.9)

\put(30,34.3){$v_2$}
\put(29.1,35.4){\circle*{1}}
\Line(29.1,35.4)(23.9,21)

\put(22,4){Figure 1: Graphs $P^{\star}$ and $F(a_1,a_2,a_3)$.}
\end{picture}
\end{center}

By Lemma \ref{b+21}, $F$ is isomorphism to either some $F(a_1,a_2,a_3)$ or $\overline{P^{\star}}$. If $F$ is isomorphism to some $F(a_1,a_2,a_3)$, the next claim determines the corresponding graph.
\begin{claim}\label{b+22}
If $F\cong F(a_1,a_2,a_3)$ for some nonnegative integers $a_1,a_2,a_3$, then $a_1=\omega, a_2=0,$ and $a_3=b-1-\omega$ satisfying $\omega \leqslant \frac{1}{2}(b-1)$.
\end{claim}
\begin{proof}
By the definition of $F(a_1,a_2,a_3)$, we have $|V_i|=a_i\ (i=1,2,3)$ (see Figure 1). If $a_2\neq 0$, let $v\in V_2$. It is easy to get that $F$ contains a double star with central vertices $v_2,v$ and $b$ leaves, a contradiction. Thus, $a_2=0$. Assume, without loss of generality, that $a_3\geqslant a_1$. If $a_1<a_3$, then $d(v_1)=a_1+1<a_3+1=d(v_3)$. By Lemma \ref{Mm}, we have $x_{v_1}<x_{v_3}$. If $a_1=a_3$, by symmetry, we have $x_{v_1}=x_{v_3}$. Hence, $x_{v_1}\leqslant x_{v_3}$. Since $\Delta(F)<b$, we have $a_1\geqslant 1$. By contracting the edge $v_1v_2$, the resulting graph contains a copy of $K_{a_1+1,a_3+1}$. Note that $F$ has the $(a,b)$-property. Then we have $a_1+1>\omega$. Thus, $\omega\leqslant a_1\leqslant a_3$. We next to show that $a_1\leqslant \omega$. Suppose that $a_1>\omega$. Let $v\in V_1$ and $F'=F-v_1v+v_3v\cong F(a_1-1,0,a_3+1)$. Since $a_1-1\geqslant\omega$, we see that $F'$ has the $(a,b)$-property. Together with Lemma \ref{ab}, $G'=G^*-E(F)+E(F')$ is $K_{a,b}$-minor free. By Lemma \ref{trans1} and $x_{v_1}\leqslant x_{v_3}$, $\lambda_\alpha(G')>\lambda(G^*)$, a contradiction. So, $a_1=\omega$. Therefore, $F\cong F(\omega,0,b-1-\omega)$ with $\omega \leqslant \frac{1}{2}(b-1)$, as desired.
\end{proof}
We come back to show Lemma~\ref{b+2}. Recall $\tau\leqslant 2$. If $\tau =1$, by $\tau=\lfloor\frac{b+1}{a+1}\rfloor$, we have $a \leqslant b \leqslant 2a$. Thus,
$$
\omega=\min \left\{a,\left\lfloor\frac{b+1}{2}\right\rfloor\right\}=\left\lfloor\frac{b+1}{2}\right\rfloor.
$$
By Claim~\ref{b+22}, suppose that $F\cong F(\omega,0,b-1-\omega)$ with $\omega\leqslant \frac{1}{2}(b-1)$. Then we have $b-1\geqslant 2\omega\geqslant b$, a contradiction. Therefore, $F\cong \overline{P^{\star}}$. If $\tau =2$, by $\tau=\lfloor\frac{b+1}{a+1}\rfloor$, we have $2a+1 \leqslant b \leqslant 3a+1$. Then,
$$
\omega=\min \left\{a,\left\lfloor\frac{b+1}{2}\right\rfloor\right\}=a.
$$
To get $F\cong F(\omega, 0, b-1-\omega)=F(a,0,b-1-a)$, we only need to show that $F\ncong \overline{P^{\star}}$. Suppose that $\overline{F}\cong P^{\star}$ (see Figure 1). By symmetry, one has $x_{v_1}=\cdots =x_{v_5}$ and $x_{u_1}=\cdots =x_{u_5}$. Note that $|V(F)|=b+2=10$ and $2a+1 \leqslant b \leqslant 3a+1$. Then we have $b=8, a=3$ and so $e(F)=30$.

By Lemma \ref{abgraph}, it is easy to get that $F(\omega, 0, b-1-\omega)=F(3,0,4)$ has the $(a,b)$-property. Note that $|V(F(3,0,4))|=10=b+2$ and $e(F(3,0,4))={7\choose 2}+7+2=30=e(F)$. Let $G'=G^*-E(\overline{P^{\star}})+E(F(3,0,4))$. By Lemma \ref{ab}, we have $G'$ is $K_{a,b}$-minor free. If $x_{u_1}>x_{v_1}$, then
$$
F(3,0,4)=\overline{P^{\star}}-\sum\limits_{i=1,i\neq2}^{5}v_2u_i-v_5u_1-v_5u_2-v_4u_3
+\sum\limits_{1\leqslant i\neq j\leqslant 5}u_iu_j+v_1u_1+v_3u_3.
$$
By some calculations, 
\begin{align*}
\lambda_\alpha(G')-\lambda_\alpha(G^*)\geqslant&
{\bf x}^{T}(A_{\alpha}(G')-A_{\alpha}(G^*)) {\bf x}\\
=&\alpha(12x_{u_1}^2+2x_{v_1}^2)+2(1-\alpha)(5x_{u_1}^2+2x_{u_1}x_{v_1})\\
&-\alpha(7x_{u_1}^2+7x_{v_1}^2)+2(1-\alpha)(7x_{u_1}x_{v_1})\\
=&5(x_{u_1}-x_{v_1})\left(\alpha x_{v_1}+(2-\alpha)x_{u_1}\right)\\
>&0,
\end{align*}
which contradicts the maximality of $\lambda_\alpha(G^*)$.

If $x_{v_1}\geqslant x_{u_1}$, via a very similar approach to that of $x_{u_1}>x_{v_1}$, we have
$$
F(3,0,4)=\overline{P^{\star}}-\sum\limits_{i=1,i\neq2}^{5}u_2v_i-u_1v_2-u_1v_4-u_3v_5+\sum\limits_{1\leqslant i\neq j\leqslant 5}v_iv_j+v_4u_4+v_5u_5,
$$
and $\lambda_\alpha(G')-\lambda_\alpha(G^*)\geqslant5(x_{v_1}-x_{u_1})\left(\alpha x_{u_1}+(2-\alpha)x_{v_1}\right)\geqslant 0.$ Suppose that $\lambda_\alpha(G')=\lambda_\alpha(G^*)$. Then we have ${\bf x}$ is the eigenvector of $G'$ and $x_{v_5}=\cdots=x_{v_1}=x_{u_1}=\cdots=x_{u_5}$. By Lemma~\ref{Mm} and $d_{F(3,0,4)}(v_1)=7>4=d_{F(3,0,4)}(u_1)$, we have $x_{v_1}>x_{u_1}$, a contradiction. Therefore, $\lambda_\alpha(G')>\lambda_\alpha(G^*)$, we get a contradiction.
\end{proof}
\begin{lem}\label{b+1}
If $F\in \mathscr{F}_{b+1}$ and $b\geqslant 4$, then $F\cong \overline{F_{a,b}}$ and $e(F)={b\choose 2}+\tau-1$ with $\tau\geqslant 2$.
\end{lem}
\begin{proof}
Note that $F$ is an edge maximal graph with the $(a,b)$-property (see Lemma \ref{emax}). According to Lemmas \ref{eb+1} and \ref{cb+1}, $e(F)={b\choose 2}+\tau-1$ and $\overline{F}$ form a forest consisting of $\tau$ components, where each component has a minimum of $\omega+1$ vertices. Our objective is to demonstrate the validity of $\overline{F}\cong F_{a,b}.$ Recall $F_{a,b}$ is a disjoint union of $\tau$ stars,  with the exception of at most one star, all of which are isomorphic to $K_{1,a}$.

Initially, we demonstrate that every element of $\overline{F}$ is a star. Assume that there is a component $F_1$ in $\overline{F}$ that is not a star. We can find two leaves in $F_1$, say $v_1,\ v_2$, such that $N_{F_1}(v_1)\cap N_{F_1}(v_2)=\emptyset$. Assume that $N_{F_1}(v_i)=\{u_i\} \ (i=1,2)$. Let $F'=F-v_1u_2+v_1u_1$ if $x_{u_1}\geqslant x_{u_2}$ and $F'=F-v_2u_1+v_2u_2$ otherwise. By Lemma \ref{cb+1}, we see that $F'$ also has the $(a,b)$-property. Let $G'=G^*-E(F)+E(F')$. By Lemmas \ref{trans1} and \ref{ab}, one sees that $G'$ is $K_{a,b}$-minor free and $\lambda_\alpha(G')>\lambda_\alpha(G^*)$, a contradiction. Thus, each component of $\overline{F}$ is a star.

Since $F$ is connected, $\overline{F}$ has at least two components and so $\left\lfloor\frac{b+1}{a+1}\right\rfloor=\tau \geqslant 2$. Then we have $\omega=\min\{a, \left\lfloor\frac{b-1}{2}\right\rfloor\}=a$. According to Lemma \ref{cb+1}, every component of $\overline{F}$ has a minimum of $a+1$ vertices. In order to get $F \cong\overline{F_{a,b}}$, it is only necessary to demonstrate that $\overline{F}$ contains at most one component that is not isomorphic to $K_{1,a}$. Otherwise, there are two components, $F_1$ and $F_2$, which satisfy $|V(F_2)|\geqslant|V(F_1)|\geqslant a+2$. Consider $v_i$ as a leaf and $u_i$ as the center vertex of $F_i\ (i=1,2)$. Then we have $d_F(u_2)\geqslant d_F(u_1)$. If $d_F(u_2)> d_F(u_1)$, by Lemma \ref{Mm}, we have $x_{u_2}>x_{u_1}$. If $d_F(u_2)=d_F(u_1)$, by symmetry, we have $x_{u_2}=x_{u_1}$. Thus, $x_{u_2}\geqslant x_{u_1}$. Let $F'=F-v_2u_1+v_2u_2$ and $G'=G^*-E(F)+E(F')$. From Lemma \ref{cb+1}, $F'$ has the $(a,b)$-property. Together with Lemma \ref{ab}, one sees that $G'$ is $K_{a,b}$-minor free. Moreover, by Lemma \ref{trans1} and $x_{u_2}\geqslant x_{u_1}$, we have $\lambda(G')>\lambda(G^*)$, a contradiction.
\end{proof}
\unitlength 1.7mm
\begin{center}
\begin{picture}(70.6,51)\linethickness{0.8pt}
\put(26,38){$K_\tau$}
\cbezier(.1,41.2)(2.4,34.6)(68.3,34.6)(70.6,41.2)\cbezier(.1,41.2)(2.4,47.7)(68.3,47.7)(70.6,41.2)

\put(5,42){$u_1$}
\put(5.3,41.1){\circle*{.8}}
\cbezier(5.3,41.1)(17.8,42.6)(30.4,43.3)(43,40.6)
\cbezier(5.3,41.1)(26,45.1)(43.2,44.3)(64.5,40.8)
\Line(5.3,41.1)(37,24.8)
\Line(5.3,41.1)(37.5,25.9)
\Line(5.3,41.1)(57.1,26.2)
\Line(5.3,41.1)(57.7,27.2)

\put(18.5,41){\circle*{.4}}
\put(24.2,41){\circle*{.4}}
\put(31.2,40.9){\circle*{.4}}

\put(43,41.5){$u_1$}
\put(43,40.6){\circle*{.8}}
\Line(43,40.6)(64.5,40.8)
\Line(43,40.6)(59.3,27.8)
\Line(43,40.6)(60.6,28.2)
\Line(43,40.6)(10,28.7)
\Line(43,40.6)(11.7,28)

\put(65,40.7){$u_{\tau}$}
\put(64.5,40.8){\circle*{.8}}
\Line(64.5,40.8)(13.3,26)
\Line(64.5,40.8)(13.9,24.5)
\Line(64.5,40.8)(43.7,28.6)
\Line(64.5,40.8)(45.3,28)

\put(3,13){$K_{b+1-\tau}$}
\polygon(0,11.8)(70.5,11.8)(70.5,30.6)(0,30.6)
\cbezier(2.5,22.2)(2.9,12.3)(15.2,12.3)(15.6,22.2)\cbezier(2.5,22.2)(2.9,32.2)(15.2,32.2)(15.6,22.2)

\put(6,21.8){$K_a$}
\put(9.5,27){$u$}
\put(10,28.7){\circle*{.8}}
\put(11.7,28){\circle*{.8}}
\put(13.3,26){\circle*{.8}}
\put(13.9,24.5){\circle*{.8}}
\cbezier(14.2,22.5)(20.9,20.1)(27.7,20.6)(37.1,21.8)
\put(14.2,22.5){\circle*{.8}}
\cbezier(13.4,21.2)(21.3,18.8)(28.1,19.3)(37.5,20.3)
\put(13.4,21.2){\circle*{.8}}
\cbezier(12.1,18.5)(26.8,13.5)(43.7,11.1)(57.7,19.5)
\put(12.1,18.5){\circle*{.8}}
\cbezier(11,17.3)(26.4,12)(43.3,10.2)(58.4,17.8)
\put(11,17.3){\circle*{.8}}

\put(18.5,23){\circle*{.4}}
\put(24.1,23){\circle*{.4}}
\put(31.3,22.9){\circle*{.4}}

\cbezier(35.8,22.2)(36.3,12.3)(48.5,12.3)(48.9,22.2)\cbezier(35.8,22.2)(36.3,32.2)(48.5,32.2)(48.9,22.2)
\put(40,21){$K_a$}
\put(37,24.8){\circle*{.8}}
\put(37.5,25.9){\circle*{.8}}
\put(43.7,28.6){\circle*{.8}}
\put(45,26){$u$}
\put(45.3,28){\circle*{.8}}
\put(37.1,21.8){\circle*{.8}}
\put(37.5,20.3){\circle*{.8}}
\Line(47.8,22.8)(55.5,22.8)
\put(47.8,22.8){\circle*{.8}}
\Line(47.9,21.5)(55.8,21.5)
\put(47.9,21.5){\circle*{.8}}

\cbezier(54.2,21.8)(54.6,11.8)(66.8,11.8)(67.3,21.8)\cbezier(54.2,21.8)(54.6,31.8)(66.8,31.8)(67.3,21.8)
\put(60,20.4){$K_c$}
\put(57.1,26.2){\circle*{.8}}
\put(57.7,27.2){\circle*{.8}}
\put(59.3,27.8){\circle*{.8}}
\put(60.6,28.2){\circle*{.8}}
\put(55.5,22.8){\circle*{.8}}
\put(55.8,21.5){\circle*{.8}}
\put(57.7,19.5){\circle*{.8}}
\put(58.4,17.8){\circle*{.8}}

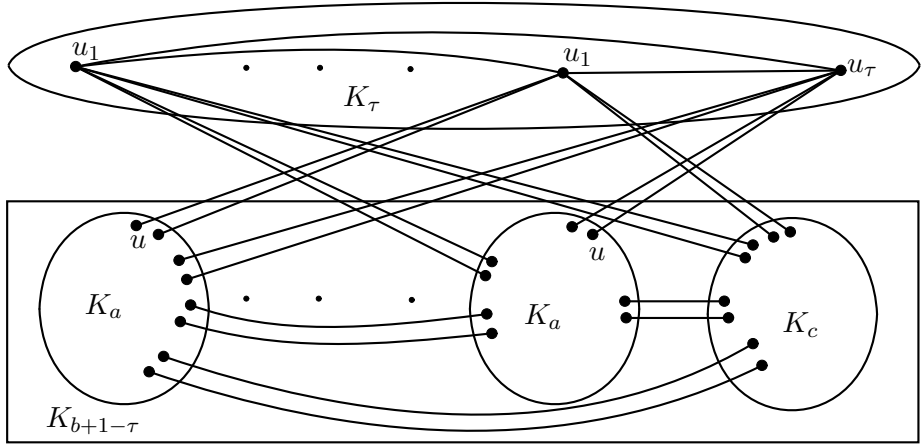
\put(12,4){Figure 2: Graph $\overline{F_{a,b}}$ for $\tau=\lfloor\frac{b+1}{a+1}\rfloor$.}
\end{picture}
\end{center}
\begin{proof}[\bf The proof of Theorem \ref{main2}] By Lemma~\ref{S}, $G^*$ has a clique dominating set $S^*$ of cardinality $a-1$. For $G^*-S^*$, by $b\geqslant a\geqslant 2$, we consider the following cases.

{\bf Case 1.}\ $2\leqslant b\leqslant 3$. In this case, if $b=2$, by Lemma \ref{b23}, we have $|\mathscr{F}_{\neq 2}|=|\mathscr{F}_1|\leqslant 1$. Thus, $G^*-S^*\cong kK_2\cup K_t \ (0\leqslant t\leqslant 1)$. If $b=3$, by Lemma \ref{b23}, we have $|\mathscr{F}_{\neq 3}|\leqslant 1$. To get $G^*-S^*\cong kK_2\cup K_t\ (0\leqslant t\leqslant 2)$, it is only necessary to demonstrate that $\mathscr{F}_{>3}=\emptyset$. Suppose that $\mathscr{F}_{>3}\neq\emptyset$. Let $F\in\mathscr{F}_{>3}$. By the proof of Lemma \ref{b23}, $F\cong P_{|V(F)|}=v_1v_2\ldots v_{|V(F)|}$ where $|V(F)|\geqslant 4$.

{\bf Case 1.1.}\ $|V(F)|=4$. Let $F'=F-v_1v_2+v_4v_2$ and $G'=G^*-E(F)+E(F')$. Notice that $F'\cong K_1\cup K_3$ has the $(a,b)$-property. By Lemma \ref{ab}, we have $G'$ is $K_{a,3}$-minor free. By symmetry, one has $x_{v_1}=x_{v_4}$. Together with Lemma \ref{trans1}, we have $\lambda_\alpha(G')>\lambda_\alpha(G^*)$, a contradiction.

{\bf Case 1.2.}\ $|V(F)|\geqslant 5$ is odd. Assume $|V(F)|=2l+1$. Then let $F'=F-v_{l-1}v_{l}-v_{l+2}v_{l+3}+v_{l-1}v_{l+3}+v_lv_{l+2}$ and $G'=G^*-E(F)+E(F')$. Clearly, $F'=K_3\cup P_{2l-2}$ has the $(a,b)$-property. By Lemma \ref{ab}, $G'$ is $K_{a,3}$-minor free. Moreover, by symmetry, one has $x_{l-1}=x_{l+3}$ and $x_l=x_{l+2}$. Then
\begin{align*}
\lambda_\alpha(G')-\lambda_\alpha(G^*)
\geqslant& {\bf x}^T(A_\alpha(G')-A_\alpha(G^*)){\bf x}\\
=&\alpha(x_{v_{l-1}}^2+x_{v_{l+3}}^2+x_{v_l}^2+x_{v_{l+2}}^2)+2(1-\alpha)(x_{v_{l-1}}x_{v_{l+3}}+x_{v_l}x_{v_{l+2}})\\
&-\alpha(x_{v_{l-1}}^2+x_{v_l}^2+x_{v_{l+2}}^2+x_{v_{l+3}}^2)-2(1-\alpha)(x_{v_{l-1}}x_{v_{l}}+x_{v_{l+2}}x_{v_{l+3}})\\
=&2(1-\alpha)(x_{v_{l-1}}-x_{v_{l}})^2\\
\geqslant&0.
\end{align*}
Suppose that $\lambda_\alpha(G')=\lambda_\alpha(G^*)$. By $1-\alpha>0$, $G'$ has ${\bf x}$ as its eigenvector and $x_{v_{l-1}}=x_{v_{l}}$. By symmetry of $G'$, one has $x_{v_l}=x_{v_{l+1}}=x_{v_{l+2}}$. Note that
\begin{align*}
\lambda_\alpha(G^*)x_{v_{i}}=\alpha(a-1)x_{v_i}+(1-\alpha)X_s+2\alpha x_{v_i}+(1-\alpha)(x_{v_{i-1}}+x_{v_{i+1}}),\ 2\leqslant i\leqslant 2l.
\end{align*}
We have $x_{v_1}=x_{v_2}=\cdots=x_{2l+1}$. By $x_{v_i}>0\ (2\leqslant i\leqslant 2l)$, we get
\begin{align*}
\lambda_\alpha(G^*)x_{v_1}=&\alpha(a-1)x_{v_1}+(1-\alpha)X_{s}+\alpha x_{v_1}+(1-\alpha)x_{v_2}\\
\neq&\alpha(a-1)x_{v_2}+(1-\alpha)X_{s}+2\alpha x_{v_2}+(1-\alpha)(x_{v_1}+x_{v_3})
=\lambda_\alpha(G^*)x_{v_2}.
\end{align*}
By $\lambda_\alpha(G^*)>0$, we have $x_{v_1}\neq x_{v_2}$, a contradiction. Thus, $\lambda_\alpha(G')\neq\lambda_\alpha(G^*)$. Then $\lambda_\alpha(G')>\lambda_\alpha(G^*)$, a contradiction.

{\bf Case 1.3.}\ $|V(F)|\geqslant 5$ is even. Assume $|V(F)|=2l$. Then let $F'=F-v_{l-1}v_{l}-v_{l+2}v_{l+3}+v_{l-1}v_{l+3}+v_lv_{l+2}=K_3\cup P_{2l-3}$ and $G'=G^*-E(F)+E(F')$. By a similar discussion as that of odd $|V(F)|$, we get $\lambda_\alpha(G')>\lambda_\alpha(G^*)$ and $G'$ is $K_{a,3}$-minor free, a contradiction.

{\bf Case 2.}\ $b\geqslant 4$. By Lemma \ref{>b+2}, we have $\mathscr{F}_{>b+2}=\emptyset$. Hence, we consider $\mathscr{F}_{\leqslant b+2}$. If $\tau=1$, by Lemmas \ref{b4>b} and \ref{b+1}, we have $\mathscr{F}_{b+1}=\emptyset$ and there all but at most one component of $G^*-S^*$ are isomorphism to $K_b$. Assume that there exists a component $F\subseteq G^*-S^*$ such that $F\ncong K_b$. By $|V(G^*-S^*)|=kb+t\ (0\leqslant t\leqslant b-1)$, we have $|V(F)|=t$ or $|V(F)|=b+t$. If $t\neq 2$, then $F\cong K_t$ and we have $G^*-S^*\cong kK_b\cup K_t$. If $t=2$,  by $\tau=1$, we have $F\cong K_2$ or $\overline{P^{\star}}$. For $b\neq8$, we have $G^*-S^*\cong kK_b\cup K_2$. For $b=8$, if $F\cong K_2$, then $e(K_8\cup K_2)=29<30=e(\overline{P^{\star}})$, a contradiction. Thus, $F\cong \overline{P^{\star}}$ and we have $G^*-S^*\cong (k-1)K_b\cup \overline{P^{\star}}$.

If $\tau \geqslant 2$, to finish the proof, we need the following claims.
\begin{claim}\label{sizeb+1}
 $|\mathscr{F}_{b+1}|\leqslant 2(\tau-1)$ for $\tau\geqslant 2$.
\end{claim}
\begin{proof}
If $\mathscr{F}_{b+1}= \emptyset$, we get $|\mathscr{F}_{b+1}|=0\leqslant 2(\tau-1)$.
If $\mathscr{F}_{b+1}\neq \emptyset$, let $F\in \mathscr{F}_{b+1}$. By Lemma \ref{b+1}, we have  $F\cong \overline{F_{a,b}}$ and $e(F)={b\choose 2}+\tau-1$.

Suppose that $|\mathscr{F}_{b+1}|\geqslant 2\tau$. Let $\mathcal{F}=2\tau \overline{F_{a,b}}$ and $\mathcal{F}'=2\tau K_b\cup K_{2\tau}$. By $a\geqslant 2$ and $b\geqslant 4$, we have $2\tau\leqslant \frac{2b+2}{a+1}\leqslant \frac{2}{3}b+\frac{2}{3}<b$. Hence, $\mathcal{F}'$ has the $(a,b)$-property. Together with Lemma \ref{emax} and $|V(\mathcal{F}')|=|V(\mathcal{F})|$, we have $e(\mathcal{F}')\leqslant e(\mathcal{F})$. On the other hand, by a direct calculation, 
$$
e(\mathcal{F}')=2\tau{b\choose 2}+{2\tau\choose 2}>2\tau\left({b\choose 2}+\tau-1\right)=e(\mathcal{F}),
$$
we get a contradiction. Then $|\mathscr{F}_{b+1}|\leqslant 2\tau-1$.

Recall $F(a,b)=K_{1,c}\cup (\tau-1)K_{1,a}$, where $c=b-(a+1)(\tau-1)$. By symmetry, let $u_1$ and $u_\tau$ denote the central vertex of $K_{1,a}$ and $K_{1,c}$, respectively. Moreover, let $u$ denote non-central vertex of $K_{1,a}$. Then we get $\overline{F_{a,b}}$, see Figure 2.

To get $|\mathscr{F}_{b+1}|\leqslant 2(\tau-1)$, we just need to show that $|\mathscr{F}_{b+1}|\neq 2\tau-1$. Otherwise, let $\mathcal{F}=(2\tau-1) \overline{F_{a,b}}$. By Lemma \ref{Mm} and $d_{\mathcal{F}}(u)=b-1>b-a=d_{\mathcal{F}}(u_1)$, we have $x_u>x_{u_1}$. Recall $c=b-(a+1)(\tau-1)\geqslant a$. If $c=a$, we have $x_{u_1}=x_{u_\tau}$. If $c>a$, then $d_{\mathcal{F}}(u_\tau)=b-c<b-a=d_{\mathcal{F}}(u_1)$. By Lemma \ref{Mm}, we have $x_{u_1}>x_{u_\tau}$. Thus, $x_u>x_{u_1}\geqslant x_{u_\tau}$.

Let $\mathcal{F}'=(2\tau-1) K_b\cup K_{2\tau-1}$. Clearly, $e(\mathcal{F}')=e(\mathcal{F})$ and $|V(\mathcal{F}')|=V(\mathcal{F})|$. Without loss of generality, assume that $V(\mathcal{F}')=V(\mathcal{F})$. It is easy to check that
\begin{align}\label{xy1}
\mathcal{F}'=\mathcal{F}
-(2\tau-1)a(\tau-1)uu_\tau-(2\tau-1)(\tau-1)u_1u_\tau
+(2\tau-1)a(\tau-1)uu_1+{2\tau-1\choose 2}u_\tau u_\tau.
\end{align}
Let $G'=G^*-E(\mathcal{F})+E(\mathcal{F}')$. As $\mathcal{F}'$ has the $(a,b)$-property, together with Lemma \ref{ab}, we have $G'$ is $K_{a,b}$-minor free. Denote by ${\bf y}$ the Perron vector of $G'$. For any vertex $w\neq u_{\tau}$, by $2\tau<b$, we have $d_{\mathcal{F}'}(w)=b-1>2\tau-2=d_{\mathcal{F}'}(u_{\tau})$. Together with Lemma \ref{Mm}, we have $y_w>y_{u_\tau}$.
By Lemma \ref{xy} and \eqref{xy1}, 
\begin{align*}
{\bf x}^{T} {\bf y}(\lambda_\alpha(G')-\lambda_\alpha(G^*))
=&{\bf x}^{T}(A_{\alpha}(G')-A_{\alpha}(G^*)) {\bf y}\\
=&(2\tau-1)(\tau-1)\left(\alpha(x_{u_\tau} y_{u_\tau}+x_{u_\tau} y_{u_\tau})+(1-\alpha)(x_{u_\tau} y_{u_\tau}+x_{u_\tau} y_{u_\tau})\right)\\
&+(2\tau-1)a(\tau-1)\left((\alpha(x_uy_w+x_{u_1} y_w)+(1-\alpha)(x_uy_w+x_{u_1}y_w)\right)\\
&-(2\tau-1)(\tau-1)\left((\alpha(x_{u_1}y_w+x_{u_\tau} y_{u_\tau})+(1-\alpha)(x_{u_1}y_{u_\tau}+x_{u_\tau} y_w)\right)\\
&-(2\tau-1)a(\tau-1)\left((\alpha(x_uy_w+x_{u_\tau} y_{u_\tau})+(1-\alpha)(x_uy_{u_\tau}+x_{u_\tau} y_w)\right)\\
=&(2\tau-1)(\tau-1)\big(\alpha(a-1)(x_{u_1}y_w-x_{u_\tau} y_{u_\tau})
+(1-\alpha)((ax_u-x_{u_\tau})\\
&\times(y_w-y_{u_\tau})+(x_{u_1}-x_{u_\tau})(ay_w-y_{u_\tau}))\big)\\
>&0.
\end{align*}
Then we get $\lambda_\alpha(G')>\lambda_\alpha(G^*)$, a contradiction.
\end{proof}
\begin{claim}\label{<bb+2}
$\mathscr{F}_{<b}\cup\mathscr{F}_{b+2}=\emptyset$ if $\mathscr{F}_{b+1}\neq \emptyset$.
\end{claim}
\begin{proof}
Let $F_1\in \mathscr{F}_{b+1}$. By Lemma \ref{b+1}, we have $\tau\geqslant 2$ and $F_1\cong \overline{F_{a,b}}$ (see Figure 2). Suppose that $\mathscr{F}_{<b}\cup\mathscr{F}_{b+2}\neq\emptyset$. By $|\mathscr{F}_{<b}\cup\mathscr{F}_{b+2}|\leqslant 1$ (see Lemma \ref{b4>b}), there exists a unique component $F_2$ in $\mathscr{F}_{<b}\cup\mathscr{F}_{b+2}$.

If $F_2\in \mathscr{F}_{<b}$, it is easy to get $F_2\cong K_{|V(F_2)|}$. Let $|V(F_2)|=h$. By Lemma~\ref{b4>b}(iii), we have $F_2\cup hK_b\subseteq G^*-S^*$. Clearly, $F_2\cup hK_b, hF_1$ have the $(a,b)$-property and $|V(F_2\cup hK_b)|=h+bh=h(b+1)=h|V(F_1)|$. By Lemma \ref{emax}, $e(hF_1)\leqslant e(F_2\cup hK_b)$, i.e.,
\[\label{VF}
h\left({b\choose 2}+\tau-1\right)\leqslant {h\choose 2}+h{b\choose 2}.
\]
By a direct calculation, we have $2\tau-1\leqslant h$. On the other hand, $F_1\cup F_2, K_b\cup K_{h+1}$ have the $(a,b)$-property and $|V(F_1\cup F_2)|=|V(K_b\cup K_{h+1})|$. Together with Lemma \ref{emax}, we have $e(F_1\cup F_2)\geqslant e(K_b\cup K_{h+1})$. If $h>\tau-1$, then
$$
e(F_1\cup F_2)={b\choose 2}+\tau-1+{h\choose 2}
<{b\choose 2}+{h\choose 2}+h
=e(K_b\cup K_{h+1}),
$$
a contradiction. So, $h\leqslant \tau-1$. By the above discussions, we have $2\tau-1\leqslant h\leqslant \tau-1$, a contradiction.

If $F_2\in \mathscr{F}_{b+2}$, by Lemma \ref{b+2}, we have $\tau=2$ and $F_2\cong F(a,0,b-1-a)$ with $b-1-a\geqslant a$ (see Figure~1). Let $\mathcal{F}=F_1\cup F_2$. Note that $d_{F_1}(u_1)=d_{F_2}(v_3)=b-a\geqslant a+1=d_{F_2}(v_1)= d_{F_1}(u_2)$. If $b-a=a+1$, by symmetry, we have $x_{u_1}=x_{u_2},\ x_{v_1}=x_{v_3}$. If $b-a>a+1$, by Lemma \ref{Mm}, we have $x_{u_1}>x_{u_2},\ x_{v_3}>x_{v_1}$. So, $x_{u_1}\geqslant x_{u_2}$ and $x_{v_3}\geqslant x_{v_1}$. For any vertex $u\in N_{F_1}(u_2)\backslash \{u_1\}$ and $v\in N_{F_2}(v_1)\backslash \{v_2\}$, we have $d_{F_1}(u)=b-1>b-a=d_{F_1}(u_1)$ and $d_{F_2}(v)=b-1>b-a=d_{F_2}(v_3)$. Together with Lemma \ref{Mm} and $d_{F_2}(v_1)=a+1>2=d_{F_2}(v_2)$, we have $x_u>x_{u_1}\geqslant x_{u_2}$ and $x_v>x_{v_3}\geqslant x_{v_1}>x_{v_2}$.

Let $\mathcal{F'}=2K_b\cup K_3$. We have $e(\mathcal{F})=e(\mathcal{F'})$ and $|V(\mathcal{F})|=|V(\mathcal{F'})|$. Without loss of generality, assume that $V(\mathcal{F})=V(\mathcal{F'})$. It is easy to check that
\begin{align}\label{xy2}
\mathcal{F'}=\mathcal{F}-u_1u_2-v_3v_2+u_2v_2+u_2v_1-a\cdot vv_1+ a\cdot vv_3-a\cdot uu_2+a\cdot uu_1.
\end{align}
Let $G'=G^*-E(\mathcal{F})+E(\mathcal{F'})$. As $\mathcal{F}'$ has the $(a,b)$-property, together with Lemma \ref{ab}, we have $G'$ is $K_{a,b}$-minor free. Denote by ${\bf y}$ the Perron vector of $G'$. For the vertex of $K_b$ and $K_3$, we use $y_1$ and $y_2$ to represent its corresponding coordinate of ${\bf y}$, respectively. Combing with Lemma \ref{Mm} and $b\geqslant 4$ gives us $y_1>y_2$.

By Lemma \ref{xy} and \eqref{xy2}, 
\begin{align*}
{\bf x}^{T} {\bf y}(\lambda_\alpha(G')-\lambda_\alpha(G^*))
=&{\bf x}^{T}(A_{\alpha}(G')-A_{\alpha}(G^*)) {\bf y}\\
=&\alpha(x_{u_2}y_2+x_{v_2}y_2+x_{u_2}y_2+x_{v_1}y_2)
+(1-\alpha)(x_{u_2}y_2+x_{v_2}y_2+x_{u_2}y_2\\
&+x_{v_1}y_2)-\alpha(x_{u_1}y_1+x_{u_2}y_2+x_{v_3}y_1+x_{v_2}y_2)
-(1-\alpha)(x_{u_1}y_2+x_{u_2}y_1\\
&+x_{v_3}y_2+x_{v_2}y_1)+a\big[
 \alpha(x_vy_1+x_{v_3}y_1)+(1-\alpha)(x_vy_1+x_{v_3}y_1)\\
&-\alpha(x_vy_1+x_{v_1}y_2)-(1-\alpha)(x_vy_2+x_{v_1}y_1)
+\alpha(x_uy_1+x_{u_1}y_1)\\
&+(1-\alpha)(x_uy_1+x_{u_1}y_1)-\alpha(x_uy_1+x_{u_2}y_2)-(1-\alpha)(x_uy_2+x_{u_2}y_1)
\big]\\
=&\alpha(a-1)((x_{u_1}+x_{v_3})y_1-(x_{u_2}+x_{v_1})y_2)
+(1-\alpha)\big[(a y_1-y_2)(x_{u_1}-x_{u_2}\\
&+x_{v_3}-x_{v_1})+(y_1-y_2)(ax_u-x_{u_2}+ax_v-x_{v_2})\big]\\
>&0.
\end{align*}
So, $\lambda_\alpha(G')>\lambda_\alpha(G^*)$, a contradiction.
\end{proof}

\begin{claim}\label{<b}
$\mathscr{F}_{<b}=\emptyset$ if $t\leqslant 2(\tau-1)$.
\end{claim}
\begin{proof}
If $\mathscr{F}_{<b}\not=\emptyset$ for $t\leqslant 2(\tau-1)$, then there exists a component $F\in \mathscr{F}_{<b}$. By a similar discussion as that of \eqref{VF}, we have $|V(F)|\geqslant 2\tau-1$. Note that $\mathscr{F}_{<b}\cup \mathscr{F}_{b+2}\neq \emptyset$. Together with Lemma~\ref{b4>b} and Claim \ref{<bb+2}, one has $|\mathscr{F}_{>b+2}|=|\mathscr{F}_{b+1}|=0$ and $|\mathscr{F}_{<b}|=1$. Recall $|V(G^*-S^*)|=kb+t\ (0\leqslant t\leqslant b-1)$. Then we have $t=|V(F)|\geqslant 2\tau-1$, which contradicts $t\leqslant 2(\tau-1)$.
\end{proof}
We consider the following cases.

{\bf Case 2.1.}\ $t\neq2$. Suppose that $\mathscr{F}_{b+2}\neq\emptyset$. By Lemma~\ref{b4>b} and Claim~\ref{<bb+2}, we have $|\mathscr{F}_{>b+3}|=|\mathscr{F}_{b+1}|=|\mathscr{F}_{<b}|=0$ and $|\mathscr{F}_{b+2}|=1$. Recall $|V(G^*-S^*)|=kb+t\ (0\leqslant t\leqslant b-1)$. Then we have $t=2$, a contradiction. Therefore, $\mathscr{F}_{b+2}=\emptyset$. If $t\leqslant 2(\tau-1)$, by Claim~\ref{<b}, we have $\mathscr{F}_{<b}=\emptyset$. So, $G^*-S^*=(k-t)K_b\cup t\overline{F_{a,b}}$. If $t\geqslant 2\tau-1$, by Claim~\ref{sizeb+1}, we obtain $\mathscr{F}_{b+1}=\emptyset$. Thus, $G^*-S^*=kK_b\cup K_t$.

{\bf Case 2.2.}\ $t=2$. By $\tau\geqslant 2$, one has $t\leqslant 2(\tau-1)$. Together with Claim \ref{<b}, we have $\mathscr{F}_{<b}=\emptyset$. If $\tau>2$, by Lemma~\ref{b+2}, we have $\mathscr{F}_{b+2}=\emptyset$. Then $|\mathscr{F}_{b+1}|=t=2$ and so $G^*-S^*=kK_b\cup 2\overline{F_{a,b}}$.

If $\tau=2=t$, suppose that $\mathscr{F}_{b+1}\neq\emptyset$. Let $\mathcal{F}=\overline{F_{a,b}}\cup \overline{F_{a,b}}$ and let $u_1,u_1',u_2,u_2'\in V(\mathcal{F})$ such that $d_{\mathcal{F}}(u_1)=d_{\mathcal{F}}(u_1')=b-1-a,\ d_{\mathcal{F}}(u_2)=d_{\mathcal{F}}(u_2')=a$. By the proof of Claim \ref{<bb+2}, we have $x_u>x_{u_1}\geqslant x_{u_2}$ for all $u\in N_{\mathcal{F}}(u_2)\backslash \{u_1\}$. By symmetry, we get $x_{u_1}=x_{u_1'},\ x_{u_2}=x_{u_2'}$.

Let $\mathcal{F'}=K_b\cup F(a,0,b-1-a)$. Clearly, $e(\mathcal{F})=e(\mathcal{F'})$ and $|V(\mathcal{F})|=|V(\mathcal{F'})|$. Without loss of generality, assume that $V(\mathcal{F})=V(\mathcal{F'})$. It is easy to get
\begin{align}\label{xy3}
\mathcal{F}'=\mathcal{F}-u_1u_2-u_1'u_2'+u_1'u_2+u_2'u_2-a\cdot uu_2+a\cdot uu_1.
\end{align}
Let $G'=G^*-E(\mathcal{F})+E(\mathcal{F'})$. As $\mathcal{F}'$ has the $(a,b)$-property, together with Lemma \ref{ab}, we have $G'$ is $K_{a,b}$-minor free. Denote by ${\bf y}$ the Perron vector of $G'$. By a similar proof as that of Claim \ref{<bb+2}, we have $y_{u_1'}\geqslant y_{u_2'}>y_{u_2}$. For the vertex of $K_b$, we use $y_1$ to represent its corresponding coordinate of ${\bf y}$. Note that $d_{\mathcal{F'}}(u_1')=b-a<b-1$. Together with Lemma \ref{Mm}, one has $y_1>y_{u_1'}$. So, $y_1>y_{u_1'}\geqslant y_{u_2'}>y_{u_2}$.

By Lemma \ref{xy} and \eqref{xy3}, 
\begin{align*}
{\bf x}^{T} {\bf y}(\lambda_\alpha(G')-\lambda_\alpha(G^*))
=&\, {\bf x}^{T}(A_{\alpha}(G')-A_{\alpha}(G)) {\bf y}\\
=&\,\alpha(x_{u_1'}y_{u_1'}+x_{u_2}y_{u_2}+x_{u_2'}y_{u_2'}+x_{u_2}y_{u_2})
+(1-\alpha)(x_{u_1'}y_{u_2}+x_{u_2}y_{u_1'}\\
&+x_{u_2}y_{u_2'}+x_{u_2'}y_{v_2})-\alpha(x_{u_1}y_1+x_{u_2}y_{u_2}+x_{u_1'}y_{u_1'}+x_{u_2'}y_{u_2'})
\\
&-(1-\alpha)(x_{u_1}y_{u_2}+x_{u_2}y_1+x_{u_1'}y_{u_2'}+x_{u_2'}y_{u_1'})+a\big(
\alpha(x_uy_1+x_{u_1}y_{u_1})\\
&+(1-\alpha)(x_vy_{u_1}+x_{u_1}y_1)-\alpha(x_uy_1+x_{u_2}y_{u_2})-(1-\alpha)(x_uy_{u_2}+x_{u_2}y_1)
\big)\\
=&\,\alpha(a-1)(x_{u_1}y_1-x_{u_2}y_{u_2})
+(1-\alpha)((ax_{u}-x_{u_2})(y_1-y_{u_2})\\
&+(x_{u_1}-x_{u_2})(ay_1-y_{u_2'}))>0.
\end{align*}Thus, $\lambda_\alpha(G')>\lambda_\alpha(G^*)$, a contradiction. Then we have $\mathscr{F}_{b+1}=\emptyset$. Combing with $\mathscr{F}_{<b}=\emptyset$ and $|V(G^*-S^*)|=kb+2$ gives us $|\mathscr{F}_{b+2}|=1$. So, $G^*-S^*=(k-1)K_b\cup F(a,0,b-1-a)$.
\end{proof}
\section{\normalsize Concluding remarks}
This paper explores the characterization of $n$-vertex $K_{a,b}$-minor free graphs having the largest $A_{\alpha}$-spectral radius. Our findings enhance and broaden the implications of Theorem~\ref{Zhai-Lin-1}. Indeed, our primary findings in this paper allow us to derive several established results from previous studies.

In Theorem~\ref{main1}, let $\alpha=0$ and $a=2$, we deduce the result obtained by Nikiforov~\cite{Niki-2t}.
\begin{cor}[\cite{Niki-2t}]
Let $n-1=kb+t\, (0\leqslant t\leqslant b-1)$ and $G$ be a $K_{2,b}$-minor free graph of order $n$. If $n$ is large enough, then $\lambda_0(G) \leqslant \lambda_0(K_1 \vee (kK_b\cup K_t))$ equality holds if and only if $G \cong K_1 \vee (kK_b\cup K_t)$.
\end{cor}

In Theorem~\ref{main1}, if $\alpha=0$, then we may characterize the $K_{a,b}$-minor free graph having the largest adjacency spectral radius for $a+b=6$, which was obtained by Wang, Chen and Fang \cite{Wang}.
\begin{cor}[\cite{Wang}]
Let $n-1=kb+t\ (0\leqslant t\leqslant b-1)$, and $(a,b)=(3,3)$ or $(a,b)=(2,4)$. If $G$ is an $n$-vertex $K_{a,b}$-minor free graph with sufficient large $n$, then $\lambda_0(G) \leqslant \lambda_0(K_{a-1}\vee(kK_b\cup K_t))$, where the equality holds if and only if $G \cong K_{a-1}\vee(kK_b\cup K_t)$.
\end{cor}

If $\alpha=0$ and $2 \leqslant a \leqslant b$, then by Theorem~\ref{main1} we may deduce Theorem~\ref{Zhai-Lin-1}. In Theorem \ref{main1}, let $\alpha=\frac{1}{2}$. Then we identify all the $K_{a,b}$-minor free graphs having the largest $Q$-spectral radius for $2 \leqslant a \leqslant b$. It was obtained by Zhang and Lou \cite{Lou-2} recently.
\begin{cor}[\cite{Lou-2}]
Let $2 \leqslant a \leqslant b, \tau=\lfloor\frac{a+1}{b+1}\rfloor$ and let $G^*$ be a graph of order $n$ that does not include any $K_{a,b}$-minor and has the maximum $Q$-spectral radius, where $n-a+1=kb+t$ is large enough and $0\leqslant t\leqslant b-1$.
\begin{wst}
\item[{\rm (i)}] If $t\leqslant 2\tau-2$, then $G^* \cong K_{a-1} \vee F_1$, where $F_1\cong (k-1)K_b\cup F(a,0,b-1-a)$ if $t=\tau=2$ and $F_1\cong (k-t)K_b\cup t\overline{F_{a,b}}$ otherwise;
\item[{\rm (ii)}]If $t\geqslant 2\tau-1$, then $G^* \cong K_{a-1} \vee F_2$, where $F_2\cong (k-1)K_b\cup \overline{P^{\star}}$ if $ t=2,\tau=1, b=8$ and $F_2\cong kK_b\cup K_t$ otherwise.
\end{wst}
\end{cor}
Moreover, let $\alpha=\frac{1}{2}$ in Theorem \ref{main2}. Then we may deduce a main result obtained by Zhang and Lou \cite{Lou-1}, which identifies the connected $K_{1,b}$-minor free graphs having the maximum $Q$-spectral radius.
\begin{cor}[\cite{Lou-1}]
For $n \geqslant b \geqslant 3$, let $G^{\star}$ attain the largest $Q$-spectral radius over all $K_{1, b}$-minor free connected graphs with $n$ vertices. Then $G^{\star}\cong \overline{F_{1,b}}$ if $n=b+1$ and
$G^{\star}\cong S^{n-b}(K_b)$ otherwise.
\end{cor}
In view of Theorem \ref{main2}, it is natural to consider the following problem.
\begin{pb}
For $n \geqslant b \geqslant 4$ and $n \neq b+1$, is it possible to identify the $n$-vertex connected $K_{1, b}$-minor free graphs that maximize the $A_\alpha$-spectral radius for  $\alpha\in(0,\frac{b+1}{2})?$
\end{pb}
\section*{\normalsize Acknowledgments} The authors would like to express their sincere gratitude to the referees for their very careful reading of the paper and for insightful comments and valuable suggestions, which improved the presentation of this paper.

\section*{\normalsize Data availability}
Our manuscript has no associated date.
\section*{\normalsize Conflict of interest} 
The authors declare that there is no conflict of interests regarding this publication.

\end{document}